# HITTING TIMES FOR INDEPENDENT RANDOM WALKS ON $\mathbb{Z}^D$

By Amine Asselah and Pablo A. Ferrari

*Université de Provence and IME-USP*

We consider a system of asymmetric independent random walks on $\mathbb{Z}^d$, denoted by $\{\eta_t, t \in \mathbb{R}\}$, stationary under the product Poisson measure $\nu_\rho$ of marginal density $\rho > 0$. We fix a pattern $\mathcal{A}$, an increasing local event, and denote by $\tau$ the hitting time of $\mathcal{A}$. By using a loss network representation of our system, at small density, we obtain a coupling between the laws of $\eta_t$ conditioned on $\{\tau > t\}$ for all times $t$. When $d \geq 3$, this provides bounds on the rate of convergence of the law of $\eta_t$ conditioned on $\{\tau > t\}$ toward its limiting probability measure as $t$ tends to infinity. We also treat the case where the initial measure is *close* to $\nu_\rho$ without being product.

**1. Introduction.** We consider asymmetric independent random walks (AIRW), denoted by $\{\eta_t, t \in \mathbb{R}\}$. Informally, we first draw an initial configuration $\eta_0 : \mathbb{Z}^d \to \mathbb{N}$. For $i \in \mathbb{Z}^d$, $\eta_0(i)$ represents the number of *particles* at site $i$ at time 0. Then, independently of each other, particles perform continuous-time random walks with transition function $p(\cdot, \cdot)$ with a nonvanishing drift $\sum_i p(0, i) i \neq 0$. For each $\rho > 0$, the AIRW process is stationary under $\nu_\rho$, a product over $\mathbb{Z}^d$, of Poisson measures of marginal intensity $\rho > 0$.

It is convenient to realize the trajectories of the stationary AIRW process as a marked Poisson process obtained as follows: (i) for each $i \in \mathbb{Z}^d$, draw $N_i$ according to a Poisson variable of intensity $\rho$; (ii) mark each particle at $i$ with a trajectory $\{\gamma_t, t \in \mathbb{R}\}$ drawn from $dP_{0,i}(\gamma)$, where we write $dP_{s,i}(\gamma)$ for the law of a continuous-time random walk $\{\gamma_t, t \in \mathbb{R}\}$, of transition $p(\cdot, \cdot)$ such that $\gamma_s = i$. We denote by $\Gamma$ a realization of such a marked Poisson process, and its law, denoted by $\mathbb{P}$, is of intensity

$$dP_\rho(\gamma) = \rho \sum_{i \in \mathbb{Z}^d} dP_{0,i}(\gamma).$$

We retrieve the occupation number $\eta_t$ by taking a *time-slice* of $\Gamma$:

$$\forall i \in \mathbb{Z}^d \qquad \eta_t(\Gamma)(i) = |\{\gamma \in \Gamma : \gamma_t = i\}|,$$

and $\{\eta_t(\Gamma), t \in \mathbb{R}\}$ is a stationary process with respect to $\nu_\rho$. On the configuration space, there is a natural order: $\eta \prec \zeta$, if for all $i \in \mathbb{Z}^d$, $\eta(i) \leq \zeta(i)$. Accordingly, we say that an event $\mathcal{A}$ is increasing if $\eta \in \mathcal{A}$ and $\eta \prec \zeta$ imply that $\zeta \in \mathcal{A}$.

We consider an increasing event $\mathcal{A}$ with support in a finite subset of $\mathbb{Z}^d$ (local), and denote by $\tau$ the hitting time of $\mathcal{A}$. We are concerned with *sharp* asymptotics of the tail distribution of $\tau$. Indeed, the existence of the following limit is obtained as a simple consequence of subadditivity (see [3])

$$\lambda(\rho) = -\lim_{t \to \infty} \frac{1}{t} \log(P_{\nu_\rho}(\tau > t)), \tag{1.1}$$

where, for a probability measure $\mu$, $P_\mu$ denotes the law of $\{\eta_t, t \geq 0\}$ with $\eta_0$ drawn from $\mu$, and $E_\mu[\cdot]$ denotes the corresponding expectation. Also, when the drift is nonzero, it is proved in [3] that $\lambda(\rho) > 0$ in any dimensions, whereas when $d \geq 3$, uniform regularity estimates are established for the law of $\eta_t$ conditioned on $\{\tau > t\}$ with initial measure $\nu_\rho$, that we denote by $T_t(\nu_\rho)$. More precisely,

$$\sup_{t \geq 0} \int \left(\frac{dT_t(\nu_\rho)}{d\nu_\rho}\right)^p d\nu_\rho < \infty \qquad \text{for any integer } p. \tag{1.2}$$

Thus, if we call $\mathcal{L}_W$ the generator of the AIRW process, and assume $d \geq 3$, [3] shows that a principal eigenfunction for $\mathcal{L}_W$, denoted by $u$, with Dirichlet boundary on $\mathcal{A}$, is obtained by taking the limit of convex combinations of $\{T_t(\nu_\rho), t \geq 0\}$ [and similarly of the dual process in $L^2(\nu_\rho)$ whose generator and principal eigenvector we denote resp. by $\mathcal{L}_W^*$ and $u^*$]. Thus, $u, u^* \in L^p(\nu_\rho)$ for any integer $p$, are decreasing with the following normalization:

$$u, u^* \geq 0 \quad \text{and} \quad \int u \, d\nu_\rho = \int u^* \, d\nu_\rho = 1.$$

Note that when $d \geq 3$, $\lambda(\rho)$ given in (1.1) is the principal eigenvalue of $\mathcal{L}_W$ with Dirichlet boundary on $\mathcal{A}$, and a variational formula for $\lambda(\rho)$ is obtained in [2].

However, neither the existence of $\mu_\rho = \lim_{t \to \infty} T_t(\nu_\rho)$, the so-called *Yaglom limit*, nor sharp asymptotics for the tail distribution of $\tau$ were established. Beside this goal, we are interested in approaching such results as obtained in an irreducible context by the Birkhoff–Hopf's theorem (see [6] for a statement in a general context and a simple proof). That is, for instance, to find explicit constants $\beta$ and $M$ such that for $\varphi$ (in certain cones) with $\varphi|_\mathcal{A} \equiv 0$

$$\int \left| e^{\lambda(\rho)t} \bar{S}_t(\varphi)(\eta) - \frac{\int \varphi u^* \, d\nu_\rho}{\int uu^* \, d\nu_\rho} u(\eta) \right| d\nu_\rho(\eta) \leq M e^{-\beta t}, \tag{1.3}$$

where $\{\bar{S}_t, t \geq 0\}$ is the semi-group of the AIRW process stopped when hitting $\mathcal{A}$.

Recently, one of the authors [2] obtained for $d \geq 3$ some $L^1(\nu_\rho)$ convergence for the Cesaro mean of $\bar{S}_t(\mathbb{1}\{\mathcal{A}^c\})$ for a larger class of processes (the so-called monotone zero-range). This was based on uniqueness for the principal Dirichlet eigenfunction in a *natural* class, and an ergodic theorem. Thus, rate of convergence escaped his approach. In the present work, we are using a special representation of the AIRW process at low density to obtain explicitly computable constant $\beta$ and $M$ such that for a probability measure $\mu$ close enough to $\nu_\rho$, and for a nonnegative function $g$ with bounded *oscillation* [see (1.12) for a precise statement],

$$(1.4) \qquad \left| E_\mu[g(\eta_t)|\tau > t] - \int g u^* \, d\nu_\rho \right| \leq M e^{-\beta t} \|g\|_{\nu_\rho},$$

where $\|\cdot\|_{\nu_\rho}$ denotes the $L^2(\nu_\rho)$-norm.

To avoid unnecessary length, we have written all our results and proofs for events of the type $\mathcal{A} := \{\eta(0) > L\}$ for an integer $L$. Any local increasing event can be treated by straightforward adaptation of the arguments with more intricate notation and expressions. More importantly, to realize the stationary AIRW process as a marked Poisson process allows us to treat also increasing *space–time* patterns. Indeed, fix $\Lambda \subset \mathbb{Z}^d$ finite, and $T > 0$. Let $\{\mathcal{A}_i, i \in \mathbb{N}\}$ be increasing events with support in $\Lambda$, and $\{t_i, i \in \mathbb{N}\}$ be an increasing subdivision of $[-T, 0]$. Then, we form

$$(1.5) \qquad \mathcal{A}_T := \{\Gamma : \eta_s(\Gamma) \in \mathcal{A}_i, \ \forall s \in [t_i, t_{i+1}[, \ \forall i \in \mathbb{N}\}.$$

We define the time-shift of a trajectory $\gamma$ by $(\theta_t \gamma)_s = \gamma_{s+t}$, and denote by $\theta_t \Gamma$ the set of trajectories of $\Gamma$ shifted by $t$. Then, we define the first occurrence of $\mathcal{A}_T$ as

$$\tau = \inf\{t \geq T : \theta_t(\Gamma) \in \mathcal{A}_T\}.$$

Our results hold also for the pattern $\mathcal{A}_T$ when the density is small enough and when $d \geq 3$.

Before stating our first result, we recall (see the Appendix) that if $\Phi(z) := \sum_i p(0, i) \exp(z \cdot i)$ for $z \in \mathbb{R}^d$, then, under mild assumptions on the transition function (see Section 2.1), we have that $0 < \inf \Phi < 1$. Also, our random walk is transient, so that if $H_0$ is the first return time to 0, then $P_{0,0}(H_0 = \infty) > 0$. Finally, if $\sigma(\gamma)$ is the diameter of the set of times, the walk $\{\gamma_t, t \in \mathbb{R}\}$ spends on site $\{0\}$, and if $\beta_d = \min(1 - \inf \Phi, P_{0,0}(H_0 = \infty)) > 0$, we see in Lemma A.3 of the Appendix that, for any fixed $0 < \beta_1 < \beta_d$, we can define the following density threshold:

$$(1.6) \qquad \rho_c(\beta_1) = \frac{1}{\int [e^{\beta_1 \sigma(\gamma)}(1 + \sigma(\gamma) + \sigma(\gamma)^2)] \, dP_{0,0}(\gamma)}.$$

We now can state a result concerning the tail asymptotics of the hitting time.

PROPOSITION 1.1. *Assume that $d \geq 3$. For any $\beta_1 < \beta_d$, and any $\rho < \rho_c(\beta_1)$, there is a number $M(\beta_1)$, such that*

$$(1.7) \qquad \left| e^{\lambda(\rho)t} P_{\nu_\rho}(\tau > t) - \frac{1}{\int uu^* \, d\nu_\rho} \right| \leq M(\beta_1) \exp(-\beta_1 t).$$

Next, we consider the Yaglom limit.

PROPOSITION 1.2. *Under the same hypotheses as Proposition 1.1, the Yaglom limit $\mu_\rho$ exists. Furthermore,*

$$(1.8) \qquad \frac{dT_t(\nu_\rho)}{d\nu_\rho} \overset{t \to \infty}{\longrightarrow} u^* = \frac{d\mu_\rho}{d\nu_\rho} \qquad \text{in } L^2(\nu_\rho).$$

*The same results hold for the dual process.*

REMARK 1.3. Our approach consists in coupling $\{T_t(\nu_\rho), t \geq 0\}$ through a loss network dynamics, as developed for contour models in [7]. We actually construct one probability space on which all the conditional laws $\{T_t(\nu_\rho), t \geq 0\}$ are realized at once, as well as the limiting object $\mu_\rho$. Moreover, the convergence, in this large space, is in the almost sure sense.

Next, we consider an initial measure which is not a product Poisson measure, but is "sandwiched" between two product Poisson measures. To formulate the next result, we need more notation. We denote by $P^*_{s,i}$ the law of a walk $\{\gamma_t, t \in \mathbb{R}\}$ with *dual* transition function $\{p^*(k,j) := p(j,k), k, j \in \mathbb{Z}^d\}$ conditioned on $\gamma_s = i$.

PROPOSITION 1.4. *Assume that $d \geq 3$ and $\rho < \rho_c(\beta_1)$. Let $C_\alpha$ be a positive constant and $\{\alpha_i, i \in \mathbb{Z}^d\}$ be positive reals with*

$$(1.9) \qquad 0 \leq 1 - \alpha_i/\rho \leq C_\alpha P^*_{0,i}(H_0 < \infty) \qquad \forall i \in \mathbb{Z}^d.$$

*Define $\nu_\alpha$ to be a product Poisson measure of marginal intensity $\alpha_i$ at $i \in \mathbb{Z}^d$, that is,*

$$(1.10) \qquad \begin{array}{l} \text{for } \eta, \nu_\rho\text{-almost surely} \\[4pt] \psi_\alpha(\eta) := \displaystyle\prod_{i \in \mathbb{Z}^d} e^{(\rho - \alpha_i)} \left( \frac{\alpha_i}{\rho} \right)^{\eta(i)} \quad \text{and} \quad d\nu_\alpha := \psi_\alpha \, d\nu_\rho. \end{array}$$

*Let $\mu_G$ be a finite range Gibbs measure (see Section 2.3) with $d\mu_G := f_G \, d\nu_\rho$ and $f_G \in L^2(\nu_\rho)$. Assume that (i) $f_G$ is decreasing, and (ii) $f_G/\psi_\alpha$ is increasing. Then,*

$$(1.11) \qquad \frac{dT_t(\mu_G)}{d\nu_\rho} \overset{t \to \infty}{\longrightarrow} u^* \qquad \text{in } L^2(\nu_\rho).$$

*Furthermore, assume that $g$ is nonnegative and satisfies for some constant $C_g$*

$$(1.12) \qquad 0 \leq g(\eta) - g(\mathsf{A}_i^+ \eta) \leq C_g P_{0,i}(H_0 < \infty) g(\eta) \qquad \forall\, i \in \mathbb{Z}^d,$$

*where $\mathsf{A}_i^+$ acts on configuration $\eta$ by adding a particle at site $i \in \mathbb{Z}^d$. Then, for any $\beta_1 < \beta_d$, there is $\bar{M}(\beta_1)$ such that*

$$(1.13) \qquad \left| E_{\mu_G}[g(\eta_t) | \tau > t] - \int g u^* \, d\nu_\rho \right| \leq C_g \bar{M}(\beta_1) e^{-\beta_1 t} \|g\|_{\nu_\rho}.$$

REMARK 1.5. Note that (i) and (ii) are stronger than $\nu_\alpha \prec \mu_G \prec \nu_\rho$. Also, when $d \geq 3$, the series $\sum_i P_{0,i}^*(H_0 < \infty)^2$ is finite (which is equivalent to the well-known fact that the Green function is square summable; see, e.g., [10]). Also, $\psi_\alpha$ of (1.10), and $g$ satisfying (1.12), are in $L^2(\nu_\rho)$ (see the proof of Theorem 3(c) in [4]).

REMARK 1.6. When considering space–time event $\mathcal{A}_T$ as in (1.5), the law of $\eta_t$ conditioned on $\{\tau > t\}$, and initial measure $\nu_\rho$, that we still denote by $dT_t(\nu_\rho)$ will converge to $d\mu_\rho = u^* \, d\nu_\rho$, but $u^*$ is no more an eigenfunction of $\mathcal{L}_W^*$ with Dirichlet boundaries. Also, $\{T_t(\nu_\rho), t \geq 0\}$ is no more a semigroup, and the subadditive argument giving the limit (1.1) does not hold.

Finally, a useful by-product of the loss network representation is a comparison between hitting times for different patterns at any density. Let $\Lambda$ be a finite subset of $\mathbb{Z}^d$, and denote

$$(1.14) \qquad 0_\Lambda := \{\eta \in \mathbb{N}^{\mathbb{Z}^d} : \exists\, i \in \Lambda, \eta(i) > 0\}.$$

For a subset $\Lambda$ of $\mathbb{Z}^d$, we denote by $H_\Lambda$ the return time in $\Lambda$ for a single walk. Also, we distinguish by a hat all quantities related to $0_\Lambda$.

PROPOSITION 1.7. *Assume that $\mathcal{A}$ is a local increasing event such that $\mathcal{A} \subset 0_\Lambda$. Then, for any dimensions $d$ and any density $\rho > 0$, we have, for $t \geq 0$,*

$$(1.15) \qquad \hat{T}_t(\nu_\rho) \prec T_t(\nu_\rho) \prec \nu_\rho \quad \text{and} \quad \hat{T}_t^*(\nu_\rho) \prec T_t^*(\nu_\rho) \prec \nu_\rho.$$

*As a consequence, for any integers $i, j$ and any $t, s \geq 0$,*

$$(1.16) \qquad \begin{aligned} &\int \left(\frac{d\hat{T}_t(\nu_\rho)}{d\nu_\rho}\right)^i \left(\frac{d\hat{T}_s^*(\nu_\rho)}{d\nu_\rho}\right)^j d\nu_\rho \\ &\qquad \geq \int \left(\frac{dT_t(\nu_\rho)}{d\nu_\rho}\right)^i \left(\frac{dT_s^*(\nu_\rho)}{d\nu_\rho}\right)^j d\nu_\rho. \end{aligned}$$

This result is useful since for $0_\Lambda$ everything can be computed. In Section 3.3 we first note that the conditional laws are ordered and converge: for $t < t'$,

$$\nu_\rho \succ \hat{T}_t(\nu_\rho) \succ \hat{T}_{t'}(\nu_\rho) \stackrel{t' \to \infty}{\longrightarrow} \hat{\mu}_\rho \quad \text{and}$$

(1.17)

$$\sup_t \left\| \frac{d\hat{T}_t(\nu_\rho)}{d\nu_\rho} \right\|_{\nu_\rho} = \left\| \frac{d\hat{\mu}_\rho}{d\nu_\rho} \right\|_{\nu_\rho} \quad (< \infty \text{ when } d \geq 3).$$

Properties (1.17) hold as well for the dual process, with all notation weighed down by an "$*$". Second, using $\mathcal{P}(\lambda)$ to denote a Poisson law of intensity $\lambda$, we have that, for $t \geq 0$,

$$\hat{T}_t(\nu_\rho) = \bigotimes_{i \in \mathbb{Z}^d} \mathcal{P}(\rho P^*_{0,i}(H_\Lambda > t)),$$

$$\hat{\mu}_\rho = \bigotimes_{i \in \mathbb{Z}^d} \mathcal{P}(\rho P^*_{0,i}(H_\Lambda = \infty)),$$

(1.18)

$$\int \frac{d\hat{T}_t(\nu_\rho)}{d\nu_\rho} \frac{d\hat{T}^*_t(\nu_\rho)}{d\nu_\rho} d\nu_\rho = \exp\left( \rho \sum_{i \in \mathbb{Z}^d} P^*_{0,i}(H_\Lambda \leq t) P_{0,i}(H_\Lambda \leq t) \right),$$

$$\int \left( \frac{d\hat{T}^*_t(\nu_\rho)}{d\nu_\rho} \right)^2 d\nu_\rho = \exp\left( \rho \sum_{i \in \mathbb{Z}^d} P_{0,i}(H_\Lambda \leq t)^2 \right).$$

REMARK 1.8. In the symmetric case (i.e., when $p^* = p$), uniform $L^2(\nu_\rho)$ estimates of the densities of $\{T_t(\nu_\rho), t \geq 0\}$ imply the existence of a Yaglom limit (see Lemma 2.3 of [5]). Thus, the domination (1.15) provides a simple proof of the existence of a Yaglom limit for independent random walks in $d \geq 5$ [since $\sum P_{0,i}(H_\Lambda = \infty)^2 < \infty$ only when $d \geq 5$ in the symmetric case].

The paper is organized as follows. In Section 2 we set our notation and assumptions. In Section 3 we construct the loss network representation. We have postponed the technical proofs, that the *clans* are almost surely finite, to Section 6. We treat the event $0_\Lambda$ of (1.14) in Section 3.3, where we prove Proposition 1.7, after showing (1.17) and (1.18). We consider the conditioned nonstationary AIRW process in Section 3.4. In Section 4 we bound discrepancies between different conditional laws, basing some estimates on classical random walks estimates which we have gathered in the Appendix. Finally, in Section 5 we apply the estimates on discrepancies to obtain hitting time estimates. Actually, Section 5 could be read before Section 3, if one is willing to assume Lemmas 4.1 and 4.3, as well as Corollaries 4.2 and 4.4.

**2. Notation and assumptions.**

2.1. *The single-particle random walk.* We consider a random walk on $\mathbb{Z}^d$ with transition function $\{p(i,j)\}$ satisfying the following assumptions:

- (0) $p(i,j) \geq 0$, $p(i,i) = 0$ and $\sum_i p(0,i) = 1$.
- (i) Translation invariant: $p(i,j) = p(0, j-i)$.
- (ii) Finite range: $p(i,j) = 0$ if $\sum_{k=1}^d |i_k - j_k| > R$.
- (iii) Irreducible: for any $i$, there is $n$ such that $p^n(0,i) > 0$.
- (iv) Nonzero drift: $\sum_i i p(0,i) \neq 0$.

Note that by (0) and (i) the transition kernel $p(\cdot,\cdot)$ is doubly stochastic. Thus, we can introduce a *dual* transition kernel $\{p^*(i,j),\ i,j \in \mathbb{Z}^d\}$, with $p^*(i,j) := p(j,i)$ satisfying (0)–(iv).

A continuous trajectory $\gamma$ of the random walk is an element of $\mathcal{D}(\mathbb{R}, \mathbb{Z}^d)$, the space of cadlag step functions $\gamma : \mathbb{R} \to \mathbb{Z}^d$. $\mathcal{D}(\mathbb{R}, \mathbb{Z}^d)$ endowed with the Skorohod topology, $\mathcal{S}$, is a complete separable metric space. For a trajectory $\gamma$, let $\Sigma(\gamma) := \{t : \gamma(t) \in \{0\}\}$, let $\boldsymbol{\sigma}(\gamma)$ be the closed convex hull of $\Sigma(\gamma)$, and let $\sigma(\gamma)$ be the length of $\boldsymbol{\sigma}(\gamma)$. Since the walk is transient, $\sigma(\gamma)$ is a.s. finite. The density of the law of $\sigma(\gamma)$ is denoted by $g_\sigma$, and we show in Lemma A.3 that there is $\beta_d > 0$ such that $\int \exp(\beta_1 \sigma(\gamma)) dP_{0,0}(\gamma) < \infty$, for any $\beta_1 < \beta_d$.

2.2. *The AIRW process.* The usual description of the AIRW is in terms of the evolution of the occupation number $\eta : \mathbb{Z}^d \mapsto \mathbb{N}$. To construct the semigroup, let

$$\alpha(i) = \sum_{n=0}^\infty 2^{-n} p^n(i,0) \qquad \text{and}$$

$$\text{for } \eta, \zeta \in \mathbb{N}^{\mathbb{Z}^d} \qquad \|\eta - \zeta\| = \sum_{i \in \mathbb{Z}^d} |\eta(i) - \zeta(i)| \alpha(i).$$

The configuration space is $\Omega = \{\eta : \|\eta\| < \infty\}$, and we call $\mathbb{L}$ the space of Lipshitz functions from $(\Omega, \|\cdot\|)$ to $(\mathbb{R}, |\cdot|)$, and $\mathbb{L}_b$ the subspace of $\mathbb{L}$ consisting of bounded functions. A semi-group $\{S_t, t \geq 0\}$ can be constructed on $\mathbb{L}$ with formal generator

(2.1) $$\mathcal{L}_W \varphi(\eta) := \sum_{i,j \in \mathbb{Z}^d} p(i,j) \eta(i) (\varphi(T_j^i \eta) - \varphi(\eta)),$$

where $T_j^i \eta(k) = \eta(k)$ if $k \notin \{i,j\}$, $T_j^i \eta(i) = \eta(i) - 1$ and $T_j^i \eta(j) = \eta(j) + 1$. This has been proven in [1] (see also [11] and [13], Section 2) for the more general class of zero-range process (where particles at the same site can interact). In the independent case, a construction can be realized by attaching a trajectory to each initial particle as mentioned in the Introduction.

In [13], Section 2, $\mathcal{L}_W$ is extended to a generator, again called $\mathcal{L}_W$ for convenience, on $L^2(\nu_\rho)$ for any $\rho > 0$. It is also shown that $\mathbb{L}_b$ is a core for $\mathcal{L}_W$.

The set of invariant measures for the AIRW process has as extremal points the family $\{\nu_\rho, \rho > 0\}$ of ergodic measures. The AIRW is *monotone* (also called attractive): the partial order is preserved under the evolution.

2.3. *Gibbs measures.* We associate with each finite subset $\Lambda$ of $\mathbb{Z}^d$ a bounded map $\Phi_\Lambda : \Omega \to \mathbb{R}$ depending only on $\{\eta(i), i \in \Lambda\}$. Also, we make a finite-range assumption: if $\Lambda$ is such that $\sup\{|i-j| : i,j \in \Lambda\} > R$, then $\Phi_\Lambda \equiv 0$. With the *potential* $\{\Phi_\Lambda : \Lambda \subset \mathbb{Z}^d, \text{ finite}\}$, we associate an *energy*

$$\forall X \subset \mathbb{Z}^d \qquad H_X(\eta) = \sum_{\Lambda \cap X \neq \varnothing} \Phi_\Lambda(\eta).$$

We denote by $\mathcal{G}_\rho(\Phi)$ the set of (Gibbs) probability measures such that, for any finite $X \subset \mathbb{Z}^d$, their conditioned laws on $\eta_{X^c} := \{\eta(i), i \notin X\}$, projected on $\mathbb{N}^X$, are given by

$$\forall \eta_X \in \mathbb{N}^X \qquad \frac{\exp(-H_X(\eta_X, \eta_{X^c}))}{Z_X(\eta_{X^c})} \prod_{i \in X} d\mathbb{P}(\rho)(\eta_X(i)),$$

where $Z_X(\eta_{X^c})$ is a positive normalizing constant, and $(\eta_X, \eta_{X^c})$ is the configuration of $\Omega$ equal to $\eta_X(i)$ for $i \in X$, and to $\eta_{X^c}(i)$ for $i \notin X$. Note that, for $\mu_G \in \mathcal{G}_\rho(\Phi)$,

$$\forall i \in \mathbb{Z}^d \qquad (1 + \eta(i)) \frac{d\mathsf{A}_i^+ \mu_G}{d\mu_G}(\eta) = \rho \exp\left(\sum_{\Lambda \ni i} \Phi_\Lambda(\eta) - \Phi_\Lambda(\mathsf{A}_i^+ \eta)\right).$$

Now, hypotheses (i) and (ii) in Proposition 1.4 read for the potential

$$(2.2) \qquad 0 \leq \sum_{\Lambda \ni i} (\Phi_\Lambda(\mathsf{A}_i^+ \eta) - \Phi_\Lambda(\eta)) \leq \log\left(\frac{\rho}{\alpha_i}\right),$$

with $\alpha_i \leq \rho$ and $\sum_i (\rho - \alpha_i)^2 < \infty$ by hypotheses (1.9) of Proposition 1.4. Finally, note that the measures in $\mathcal{G}_\rho(\Phi)$ are in general not translation invariant.

### 3. The loss network representation.

*On the use of loss network.* We introduce in Section 3.1.1 a marked-Poisson point process $\Gamma_\varnothing$ with law $\mathbb{P}$. Also, for each positive number $T$, we are given a compatibility condition $\mathcal{C}_T$ on $\mathcal{G}$, the space of realizations of $\Gamma_\varnothing$. The problem that we will face in Section 4 is to compare $\mathbb{P}$ conditioned on $\mathcal{C}_T$ for different $T$'s. For this purpose, we build a birth and death process on $\mathcal{G}$ reversible with respect to $\mathbb{P}$. We actually realize the birth and death dynamics as a point process on a larger space, where we incorporate birth-time and death-time of each point, thus obtaining the *noninteracting* *rectangle process* which we denote by **C** in Section 3.1.2. The key observation is that the evolution obtained by canceling births which violate $\mathcal{C}_T$ is reversible with respect to $\mathbb{P}(\cdot|\mathcal{C}_T)$. The loss network is the *trimmed rectangle process* of Section 3.2.1 where we erase all rectangles violating $\mathcal{C}_T$. The trimming algorithm requires that the distribution of points be sparse enough. Thus, when the particle density $\rho$ is small enough, the loss network yields a coupling for all $\{\mathbb{P}(\cdot|\mathcal{C}_T), T > 0\}$.

### 3.1. *The noninteracting process.*

#### 3.1.1. *Static.*
Let $m_c$ be the counting measure on $\mathbb{Z}^d$, that is, $m_c : \mathcal{P}(\mathbb{Z}^d) \to \mathbb{N} \cup \{\infty\}$ with $m_c(\Lambda) = |\Lambda|$, the number of sites in $\Lambda \in \mathcal{P}(\mathbb{Z}^d)$, the collection of subsets of $\mathbb{Z}^d$. For any density $\rho > 0$, consider the Poisson point process on $\mathbb{Z}^d$ with intensity measure $\rho dm_c$, and denote its counting variable at site $i \in \mathbb{Z}^d$ by $N_i$. Our mark-space is $(\mathcal{D}(\mathbb{R}, \mathbb{Z}^d), \mathcal{S})$ and we consider the marked-point process with counting measure

$$p_c(\omega, \Lambda \times A) := \sum_{i \in \Lambda} \sum_{k \geq 1} \mathbb{1}\{\gamma_k^{(i)}(\omega) \in A\} \mathbb{1}\{N_i(\omega) \geq k\} \qquad \text{for } A \in \mathcal{S}, \Lambda \subset \mathbb{Z}^d,$$

where for each site $i \in \mathbb{Z}^d$, $\{\gamma_k^{(i)}, k \in \mathbb{N}\}$ are i.i.d. random walks drawn from $dP_{0,i}$, independent from $\{N_i, i \in \mathbb{Z}^d\}$. Thus, $\omega \mapsto p_c(\omega, \cdot)$ is a random measure on $\mathbb{Z}^d \times \mathcal{D}(\mathbb{R}, \mathbb{Z}^d)$ with the product Borel $\sigma$-field and deterministic intensity measure

$$\lambda(\Lambda, d\gamma) := \sum_{i \in \Lambda} \rho P_{0,i}(d\gamma).$$

We define the intensity measure [on $\mathcal{D}(\mathbb{R}, \mathbb{Z}^d)$] due to all sites of $\mathbb{Z}^d$:

(3.1) $$P_\rho(d\gamma) := \lambda(\mathbb{Z}^d, d\gamma) = \rho \sum_{i \in \mathbb{Z}^d} P_{0,i}(d\gamma).$$

We show now that $P_\rho(d\gamma)$ is space and time translation invariant.

For a trajectory $\gamma \in \mathcal{D}(\mathbb{R}, \mathbb{Z}^d)$, a time $t$ and a site $i \in \mathbb{Z}^d$, we call $\theta'_i \gamma$ the space-translation by $i$, that is, $\theta'_i \gamma = \gamma + i$. For a measure $P$ on $\mathcal{D}(\mathbb{R}, \mathbb{Z}^d)$, we call $\theta_t P(A) = P(\theta_t A)$ and $\theta'_i P(A) = P(\theta'_i A)$. First, note that $\theta_s P_{t,i} = P_{t-s,i}$. Indeed, it is enough to consider

$$A = \{\gamma \in \mathcal{D}(\mathbb{R}, \mathbb{Z}^d) : \gamma_{t_n} \in U_n, n \in \mathbb{N}\},$$

where $\{t_n\}$ is a sequence in $\mathbb{R}$, and $\{U_n\}$ in $\mathcal{P}(\mathbb{Z}^d)$. Now

$$\theta_s A = \{\gamma \in \mathcal{D}(\mathbb{R}, \mathbb{Z}^d) : \gamma_{t_n+s} \in U_n\}$$

and

$$P_{t,i}(\theta_s A) = P(\{\gamma_{t_n+s} \in U_n, n \in \mathbb{N}\}|\gamma_t = i)$$
$$= P(\{\gamma_{t_n} \in U_n, n \in \mathbb{N}\}|\gamma_{t-s} = i)$$
$$= P_{t-s,i}(A),$$

since the generator of a single walk is time-independent. Now, for any $t, t'$ and $i, j \in \mathbb{Z}^d$,

$$(3.2) \qquad P_{t,i}(d\gamma)\mathbb{1}\{\gamma_{t'} = j\} = P_{t',j}(d\gamma)\mathbb{1}\{\gamma_t = i\}.$$

Indeed, first call $Q_t(i,j)$ the probability that $\gamma_t$ is in $j$ at time $t$ given it is in $i$ at time 0 and $Q_t^*(i,j)$ the semigroup of the time-reversed walk with transition function $p^*(\cdot,\cdot)$. To see (3.2), divide the first member by $Q_{t'-t}(i,j)$ and the second one by $Q_{t-t'}^*(j,i)$ (an identical quantity) to obtain in both terms the law of a random walk with transition function $p(\cdot,\cdot)$ conditioned to be in $i$ at time $t$ and in $j$ at time $t'$. Thus,

$$(3.3) \qquad \begin{aligned} \sum_i P_{t,i}(d\gamma) &= \sum_i \sum_j \mathbb{1}\{\gamma_{t'} = j\} P_{t,i}(d\gamma) \\ &= \sum_i \sum_j \mathbb{1}\{\gamma_t = i\} P_{t,j'}(d\gamma) \\ &= \sum_j P_{t',j}(d\gamma), \end{aligned}$$

where we used (3.2) in the second equality. The time-translation invariance follows at once. The space-translation invariance is obvious by definition of $dP_\rho$.

Let $\mathcal{G}$ be the space of point measures on $\mathcal{D}(\mathbb{R}, \mathbb{Z}^d)$. Let $\Gamma_\varnothing$ be the marked-Poisson process in $\mathcal{D}(\mathbb{R}, \mathbb{Z}^d)$ with intensity measure $dP_\rho$, that is, a random variable with value in $\mathcal{G}$. The index $\varnothing$ refers to the fact that trajectories are *noninteracting*. We denote by $\mathbb{P}$ and $\mathbb{E}$, respectively, the probability and expectation induced by $\Gamma_\varnothing$. By translation invariance of $P_\rho(d\gamma)$, the law of $\Gamma_\varnothing$ is invariant by time and space translations:

$$(3.4) \qquad \theta_s \theta'_j \Gamma_\varnothing \stackrel{d}{=} \Gamma_\varnothing \qquad \text{for } (s,j) \in \mathbb{R} \times \mathbb{Z}^d.$$

We associate with the Poisson process $\Gamma_\varnothing$, and at each time $t \in \mathbb{R}$, a configuration dubbed its *time-slice* defined by

$$(3.5) \qquad \eta_t(\Gamma_\varnothing)(i) := |\{\gamma \in \Gamma_\varnothing : \gamma_t = i\}|.$$

Note that by time-translation invariance, $\eta_t(\Gamma_\varnothing)$ has the same law as $\{p_c(\{i\} \times \mathcal{D}(\mathbb{R}, \mathbb{Z}^d)), i \in \mathbb{Z}^d\}$, so that $\eta_t(\Gamma_\varnothing) = \{\eta_t(\Gamma_\varnothing)(i) : i \in \mathbb{Z}^d\}$ has law $\nu_\rho$. Moreover, by the independence of the trajectories in $\Gamma_\varnothing$, each evolving with transition $p(\cdot,\cdot)$, and the translation invariant property (3.4), $\{\eta_t(\Gamma_\varnothing), t \in \mathbb{R}\}$ is a stationary Markov process with generator $\mathcal{L}_W$ and time marginal $\nu_\rho$.

3.1.2. *Birth and death process.* We define a birth and death process on $\mathcal{G}$, whose unique reversible measure is the law of $\Gamma_\varnothing$. Following [7], we define a Poisson process on $\mathcal{D}(\mathbb{R}, \mathbb{Z}^d) \times \mathbb{R} \times \mathbb{R}^+$ of intensity $P_\rho(d\gamma)\, db\, e^{-l}\, dl$. Each point of this process is a triplet $(\gamma, b, l)$ which is associated with the *rectangle*

$$(3.6) \qquad \mathsf{R} = \gamma \times [b, b+l],$$

with *basis* $\gamma$, *birth-time* $b$, *death-time* $b+l$ and *life-epoch* $[b, b+l]$. In this case we use the notation $\gamma(\mathsf{R}) = \gamma$, $\mathrm{epoch}(\mathsf{R}) = [b, b+l]$ and $\mathrm{birth}(\mathsf{R}) = b$. The random set of rectangles induced by a realization of the Poisson process is called $\mathbf{C}$, and its law is denoted by $\mathbb{Q}$.

The process $\{\Gamma_\varnothing^\mathfrak{b} : \mathfrak{b} \in \mathbb{R}\}$ with

$$(3.7) \qquad \Gamma_\varnothing^\mathfrak{b} := \{\gamma(\mathsf{R}) : \mathsf{R} \in \mathbf{C}, \mathrm{epoch}(\mathsf{R}) \ni \mathfrak{b}\}$$

is stationary, Markov and has generator

$$(3.8) \qquad \mathcal{L}_{\mathrm{bd}} f(\Gamma) = \int P_\rho(d\gamma)[f(\Gamma \cup \{\gamma\}) - f(\Gamma)] + \sum_{\gamma \in \Gamma}[f(\Gamma \setminus \{\gamma\}) - f(\Gamma)].$$

The unique invariant (and reversible) measure for this process is $\mathbb{P}$. We omit the proofs since these facts are similar to [7], proof of Theorem 3.1, and [8], Appendix A, proof of Theorem 1. Note also that, for any $\mathfrak{b}$, the time-slices $\{\eta_t(\Gamma_\varnothing^\mathfrak{b}), t \in \mathbb{R}\}$ define an AIRW process stationary with respect to $\nu_\rho$.

3.2. *The conditioned process.* Let $\mathcal{A} := \{\eta(0) > L\}$ for any integer $L$. Fix an interval $I \subset \mathbb{R}$, and define

$$(3.9) \qquad \mathcal{A}_I := \{\Gamma \in \mathcal{G} : \exists\, s \in I \text{ with } \eta_s(\Gamma) \in \mathcal{A}\}.$$

We now need to study $\mathbb{P}$ conditioned on $\mathcal{A}_I^c$, that is,

$$(3.10) \qquad d\mathbb{P}^I(\Gamma) := \frac{\mathbb{1}\{\mathcal{A}_I^c\}(\Gamma)\, d\mathbb{P}(\Gamma)}{\mathbb{P}(\mathcal{A}_I^c)}.$$

Let $\{\Gamma_I^\mathfrak{b}, \mathfrak{b} \in \mathbb{R}\}$ be a process evolving with the same dynamics as $\{\Gamma_\varnothing^\mathfrak{b}, \mathfrak{b} \in \mathbb{R}\}$, except that jumps to $\mathcal{A}_I$ are prohibited. Since $\mathbb{P}$ is reversible for $\{\Gamma_\varnothing^\mathfrak{b}, \mathfrak{b} \in \mathbb{R}\}$, one expects that $\mathbb{P}^I$ is reversible for $\{\Gamma_I^\mathfrak{b}, \mathfrak{b} \in \mathbb{R}\}$. This is shown in the next section. The formal generator of $\{\Gamma_I^\mathfrak{b}, \mathfrak{b} \in \mathbb{R}\}$ is given by

$$(3.11) \quad \begin{aligned} \mathcal{L}_{\mathrm{bd}}^I f(\Gamma) = &\int P_\rho(d\gamma)\mathbb{1}\{\Gamma \cup \{\gamma\} \in \mathcal{A}_I^c\}[f(\Gamma \cup \{\gamma\}) - f(\Gamma)] \\ &+ \sum_{\gamma \in \Gamma}[f(\Gamma \setminus \{\gamma\}) - f(\Gamma)]. \end{aligned}$$

3.2.1. *Construction of the loss network.* For each interval $I \subset \mathbb{R}$, we say that trajectories $\gamma$ and $\gamma'$ *$I$-interact* if $\boldsymbol{\sigma}(\gamma) \cap \boldsymbol{\sigma}(\gamma') \cap I \neq \varnothing$. We say that a rectangle $\mathsf{R}'$ is an *$I$-parent* of $\mathsf{R}$ if $\gamma(\mathsf{R})$ $I$-interacts with $\gamma(\mathsf{R}')$, and the birth of $\mathsf{R}$ belong to the epoch of $\mathsf{R}'$. We call $\mathbf{A}_1^{\mathsf{R}}(I)$ the set of $I$-parents of $\mathsf{R}$. Also, we define the $n$th generation of $I$-parents and the *$I$-clan* of $\mathsf{R}$, respectively, by

$$(3.12) \quad \mathbf{A}_n^{\mathsf{R}}(I) := \bigcup_{\mathsf{R}' \in \mathbf{A}_{n-1}^{\mathsf{R}}} \mathbf{A}_1^{\mathsf{R}'}(I) \quad \text{and} \quad \mathbf{A}^{\mathsf{R}}(I) := \bigcup_{n \geq 1} \mathbf{A}_n^{\mathsf{R}}(I) \cup \{\mathsf{R}\}.$$

Note that if $J \subset I$, then $\mathbf{A}^{\mathsf{R}}(J) \subset \mathbf{A}^{\mathsf{R}}(I)$. Thus, the following result needs only to be proven for $I = \mathbb{R}$.

LEMMA 3.1. *Let $\beta_1 < \beta_d$ given in* (A.9), *and $\rho_c(\beta_1) > 0$ given in* (1.6). *For any density $\rho < \rho_c(\beta_1)$, for any interval $I \subset \mathbb{R}$, the clan of each rectangle $\mathsf{R}$ in $\mathbf{C}$ is finite $\mathbb{Q}$-almost surely.*

We prove this lemma in Section 6.1.

REMARK 3.2. We could show that, for all bounded interval $I$, and any $\rho > 0$, the $I$-clan of $\mathsf{R}$ is finite for all $\mathsf{R} \in \mathbf{C}$, $\mathbb{Q}$-almost surely. Indeed, since the interaction of trajectories $\gamma$ is through $\boldsymbol{\sigma}(\gamma)$, we have to study a loss network of intervals in a finite box.

*The $I$-trimming algorithm.* Since the clan of every rectangle is finite when $\rho < \rho_c(\beta_1)$, we can order those rectangles by birth time. Iteratively, we label each rectangle of the clan as *$I$-kept* or *$I$-deleted* in the following way. Fix a time $\mathfrak{b}$ and let $\mathsf{R}$ be a rectangle of $\mathbf{C}$ alive at $\mathfrak{b}$:

- Let $\mathsf{R}_1$ be the eldest rectangle of $\mathbf{A}^{\mathsf{R}}(I)$. If $\{\gamma(\mathsf{R}_1)\} \in \mathcal{A}_I$, then $\mathsf{R}_1$ is $I$-deleted, else it is $I$-kept.
- Assume we have $I$-labeled $C_n := \{\mathsf{R}_1, \ldots, \mathsf{R}_n\}$, the eldest $n$ rectangles of the clan. If $\{\gamma(\mathsf{R}) : \mathsf{R} \in C_n, I\text{-kept}, \text{epoch}(\mathsf{R}) \ni \text{birth}(\mathsf{R}_{n+1})\} \cup \{\gamma(\mathsf{R}_{n+1})\} \in \mathcal{A}_I$, then $\mathsf{R}_{n+1}$ is $I$-deleted, else it is $I$-kept.
- Stop the labeling once $\mathsf{R}$ is labeled.

Repeating this procedure with all rectangles alive at $\mathfrak{b}$, we obtain $\mathbf{K}(I, \mathfrak{b}) \subset \mathbf{C}$, the resulting set of $I$-kept rectangles alive at $\mathfrak{b}$. The construction of $\mathbf{K}(I, \mathfrak{b})$ does not depend on the order the clans are chosen since the labeling of a rectangle depends only on its parents.

Define, for any $\mathfrak{b} \in \mathbb{R}$,

$$(3.13) \quad \Gamma_I^{\mathfrak{b}} := \{\gamma(\mathsf{R}) : \mathsf{R} \in \mathbf{K}(I, \mathfrak{b})\}.$$

PROPOSITION 3.3. *For any interval $I \subset (-\infty, 0]$, the process $\{\Gamma_I^{\mathfrak{b}} : \mathfrak{b} \in \mathbb{R}\}$ is stationary and it is Markov with generator $\mathcal{L}_{\mathrm{bd}}^I$ defined in (3.11). The time marginal law of $\Gamma_I^{\mathfrak{b}}$ is the measure $\mathbb{P}^I$.*

PROOF. Both stationarity and the fact that the process is Markov with generator $\mathcal{L}_{\mathrm{bd}}^I$ follow from the construction. Moreover, the measure $\mathbb{P}^I$ is reversible for the process. These facts are easy to check (see [7], proof of Theorem 3.1, and [8], Appendix A, proof of Theorem 1, where details are given). □

For any interval $I$, and any $s, \mathfrak{b} \in \mathbb{R}$, we define the configuration

$$\eta_s(\Gamma_I^{\mathfrak{b}})(i) = |\{\mathsf{R} \in \Gamma_I^{\mathfrak{b}} : \gamma(\mathsf{R})_s = i\}|. \tag{3.14}$$

Since the law of $\Gamma_I^{\mathfrak{b}}$ is stationary, we often drop the superscript $\mathfrak{b}$ when no confusion is possible.

LEMMA 3.4. *For $t > 0$, we set $I = [-t, 0]$. The particle configuration $\eta_0(\Gamma_I^{\mathfrak{b}})$ defined by (3.14) has law $T_t(\nu_\rho)$.*

PROOF. Using the shorthand notation $\eta. = \{\eta_s, s \in \mathbb{R}\}$, we define $\tau(\eta.) := \inf\{s > 0 : \eta_s \in \mathcal{A}\}$, and note that

$$\{\Gamma_\varnothing \notin \mathcal{A}_{[0,t]}\} = \{\Gamma_\varnothing : \eta_s(\Gamma_\varnothing) \notin \mathcal{A}, \text{ for } s \in [0,t]\} = \{\tau(\eta.(\Gamma_\varnothing)) > t\}.$$

Now, we have seen that $\theta_t \Gamma = \Gamma$ in law. Thus, since $\theta_t \mathcal{A}_I^c = \{\tau > t\}$, we have $\mathbb{P}(\mathcal{A}_I^c) = P_{\nu_\rho}(\tau > t)$, and

$$E[f(\eta_0(\Gamma_I^{\mathfrak{b}}))] = \int f(\eta_0(\Gamma)) \frac{\mathbb{1}\{\mathcal{A}_I^c\}(\Gamma) \, d\mathbb{P}(\Gamma)}{\mathbb{P}(\mathcal{A}_I^c)}$$

$$= \int f(\eta_t) \frac{\mathbb{1}\{\tau(\eta.) > t\} \, d\mathbb{P}}{P_{\nu_\rho}(\tau > t)}$$

$$= \int f \, dT_t(\nu_\rho). \qquad \square$$

We introduce now a key object. For a realization of $\mathbf{C}$ and $\mathsf{R} \in \mathbf{C}$, we introduce the width of the clan of $\mathsf{R}$, denoted $W(\mathsf{R})$, as follows:

$$W(\mathsf{R}) := \bigcup \{\boldsymbol{\sigma}(\mathsf{R}') : \mathsf{R}' \in \mathbf{A}^\mathsf{R}((-\infty, 0])\}. \tag{3.15}$$

This width is similar to the space-width of [7] (compare with the definition of $SW$ in page 917 of [7]).

REMARK 3.5. A simple observation is that, for $t \geq 0$, if $W(\mathsf{R}) \cap (-\infty, -t] = \varnothing$, then $\mathbf{A}^\mathsf{R}((-t, 0]) = \mathbf{A}^\mathsf{R}((-\infty, 0])$.

REMARK 3.6. We call $T_\infty(\nu_\rho)$ the law of $\eta_0(\Gamma^\flat_{(-\infty,0]})$. The study of discrepancies between $T_t(\nu_\rho)$ and $T_\infty(\nu_\rho)$, in Section 4.1 and Section 5, shows that $T_\infty(\nu_\rho)$ is the Yaglom limit.

3.3. *Example*: $0_\Lambda := \{\eta : \sum_\Lambda \eta(i) > 0\}$. We recall that all quantities referring to $0_\Lambda$ have a hat. The event $0_\Lambda$ is particularly nice since when trimming $\mathbf{C}$ into $\hat{\mathbf{K}}(I,\flat)$, the trajectories touching $\Lambda$ in the time-interval $I$ are $I$-deleted. Thus, there is no need for $I$-parents, and no threshold in density. Thus, for any $\rho > 0$ and $\flat \in \mathbb{R}$, the following is well defined:

$$(3.16) \qquad \hat{\mathbf{K}}(I,\flat) = \{\mathsf{R} \in \mathbf{C} : \Sigma(\mathsf{R}) \cap I = \varnothing, \text{epoch}(\mathsf{R}) \ni \flat\},$$

where $\Sigma(\mathsf{R}) := \{t \in \mathbb{R} : \gamma_t(\mathsf{R}) \in \Lambda\}$. If $I \subset J$, then $\hat{\mathbf{K}}(I,\flat) \supset \hat{\mathbf{K}}(J,\flat)$, from where $\hat{\Gamma}^\flat_I \supset \hat{\Gamma}^\flat_J$, which implies

$$(3.17) \qquad \eta_s(\hat{\Gamma}^\flat_I) \geq \eta_s(\hat{\Gamma}^\flat_J) \qquad \text{for all } s \in \mathbb{R}.$$

Taking $t' > t > 0$, $I = [-t,0]$, $J = [-t',0]$ and $s = 0$, and using Lemma 3.4, we get

$$(3.18) \qquad \nu_\rho \succ \hat{T}_t(\nu_\rho) \succ \hat{T}_{t'}(\nu_\rho) \succ \hat{\mu}_\rho := \hat{T}_\infty(\nu_\rho).$$

Note, however, that (3.18) is not true for any increasing pattern. Also, when $I = [-t,0]$, $\hat{T}_t(\nu_\rho)$ is the law of $\eta_0(\hat{\Gamma}^\flat_I)$ by Lemma 3.4. Thus, a convenient way of obtaining $\hat{T}_t(\nu_\rho)$ is to draw at time 0, at each site $i \in \mathbb{Z}^d$, a Poisson process of intensity $\rho$ associated with trajectories $\gamma$ drawn from $P_{0,i}(d\gamma)$, but to keep only the marks which satisfy $\Sigma(\gamma) \cap I = \varnothing$. In other words, at each mark—at site $i$—we toss a coin with tail probability $P^*_{0,i}(H_\Lambda > t)$, and keep the mark if tail comes up. This yields by a classical exercise

$$\hat{T}_t(\nu_\rho) = \bigotimes_{i \in \mathbb{Z}^d} \mathcal{P}(\rho P^*_{0,i}(H_\Lambda > t)),$$

and the explicit expressions (1.18) follow easily.

PROOF OF PROPOSITION 1.7. Let $\mathcal{A}$ be a local increasing pattern with $\mathcal{A} \subset 0_\Lambda$, and set $I = [-t,0]$. By Remark 3.2, the $I$-clan is always almost-surely finite when $I$ itself is finite. It is obvious that $\hat{\mathbf{K}}(I,\flat) \supset \mathbf{K}(I,\flat)$, so that

$$(3.19) \qquad \eta_0(\hat{\Gamma}^\flat_I) \prec \eta_0(\Gamma^\flat_I) \prec \eta_0(\Gamma^\flat_\varnothing) \Longrightarrow \hat{T}_t(\nu_\rho) \prec T_t(\nu_\rho) \prec \nu_\rho.$$

A similar bound holds for the dual process. Note that when the Yaglom limits exist for both conditional laws (e.g., when $d \geq 3$ and $\rho$ small enough), then (3.19) yields that $\hat\mu_\rho \prec \mu_\rho \prec \nu_\rho$, and for the dual dynamics $\hat\mu_\rho^* \prec \mu_\rho^* \prec \nu_\rho$. In general, let us call

$$\hat{u}_t = \frac{d\hat{T}_t^*(\nu_\rho)}{d\nu_\rho}, \qquad \hat{u}_t^* = \frac{d\hat{T}_t(\nu_\rho)}{d\nu_\rho},$$

(3.20)

$$u_t = \frac{dT_t^*(\nu_\rho)}{d\nu_\rho} \quad \text{and} \quad u_t^* = \frac{dT_t(\nu_\rho)}{d\nu_\rho}.$$

For any integers $i, j, k, l$, with $j \geq 1$, we note that $\hat{u}_t^i(\hat{u}_s^*)^{j-1}u_t^k(u_s^*)^l$ is decreasing. Now, we use (3.19) to obtain

$$\int \hat{u}_t^i(\hat{u}_s^*)^j u_t^k (u_s^*)^l \, d\nu_\rho = \int \hat{u}_t^i(\hat{u}_s^*)^{j-1} u_t^k (u_s^*)^l \, d\hat{T}_s(\nu_\rho)$$
$$\geq \int \hat{u}_t^i(\hat{u}_s^*)^{j-1} u_t^k (u_s^*)^l \, dT_s(\nu_\rho).$$

Note that when $i \geq 1$,

$$\int \hat{u}_t^i(\hat{u}_s^*)^{j-1} u_t^k (u_s^*)^l \, dT_s(\nu_\rho) = \int \hat{u}_t^{i-1}(\hat{u}_s^*)^{j-1} u_t^k (u_s^*)^{l+1} \, d\hat{T}_t^*(\nu_\rho),$$

and we use that $\hat{u}_t^{i-1}(\hat{u}_s^*)^{j-1} u_t^k (u_s^*)^{l+1}$ is decreasing and (3.19) (for the dual) to obtain

$$\int \hat{u}_t^{i-1}(\hat{u}_s^*)^{j-1} u_t^k (u_s^*)^{l+1} \, d\hat{T}_t^*(\nu_\rho) \geq \int \hat{u}_t^{i-1}(\hat{u}_s^*)^{j-1} u_t^k (u_s^*)^{l+1} \, dT_t^*(\nu_\rho)$$
$$= \int \hat{u}_t^{i-1}(\hat{u}_s^*)^{j-1} u_t^{k+1} (u_s^*)^{l+1} \, d\nu_\rho.$$

Combining the last two inequalities and proceeding by induction, we obtain (1.16), and conclude the proof. $\square$

3.4. *Gibbs measure as initial conditions.* We construct a birth and death process similar to that of Section 3.1.2, but with time-slice configurations at time $-t$ drawn from $\mu_G \in \mathcal{G}_\rho(\Phi)$ (see Section 2.3). First, let us define rates of birth and death satisfying detailed balance with respect to $\mu_G$. If the configuration is $\eta$, a particle is added at site $i$ with rate $c_i(\eta)$ and dies with rate 1; in terms of occupation numbers, $\eta(i)$ grows to $\eta(i) + 1$ with rate $c_i(\eta)$, whereas $\eta(i) + 1$ decreases to $\eta(i)$ with rate $\eta(i) + 1$. Thus, at each site $i \in \mathbb{Z}^d$, we choose $c_i : \Omega \to \mathbb{R}^+$ as

(3.21) $\quad c_i(\eta) = (\eta(i) + 1)\dfrac{d\mathsf{A}_i^+ \mu_G}{d\mu_G}(\eta) \qquad \left[\text{and } \dfrac{\alpha_i}{\rho} \leq c_i \leq 1 \text{ by } (2.2)\right].$

Also, since $\mu_G$ is finite range, we have $c_i(\eta) = c_i(A_j^+ \eta)$ for all $j$ such that $|i - j| > R$, where $R$ is the range of $\mu_G$.

With each rectangle R of **C**, introduced in Section 3.1.2, we associate a random variable $U$, uniform in $[0,1]$ and independent of "everything" else. Let $\mathbf{C}_U$ be the collection of couples $(\mathsf{R},U)$ just introduced, whose law we continue to call $\mathbb{Q}$ and whose elements we continue to call rectangles for simplicity.

For pedagogical reasons, we first build a configuration whose law is $T_t(\nu_\alpha)$. For any $s \in \mathbb{R}$, we define

(3.22) $$\underline{\mathbf{C}}_s := \{(\mathsf{R},U) \in \mathbf{C}_U : U \leq \alpha_i/\rho \text{ with } i = \gamma(\mathsf{R})_s\},$$

and for an arbitrary $\mathfrak{b}$, its corresponding time-slice at $t > s$

$$\underline{\eta}_t(i) := |\{(\mathsf{R},U) \in \underline{\mathbf{C}}_s : \gamma(\mathsf{R})_t = i, \text{ and } \mathrm{epoch}(\mathsf{R}) \ni \mathfrak{b}\}|.$$

Then, $\{\underline{\eta}_t, t \geq s\}$ is Markov with generator $\mathcal{L}_W$ and initial distribution $\nu_\alpha$ at time $s$.

Henceforth, we fix $t > 0$ and set $I = [-t,0]$. We apply the trimming algorithm of Section 3.2.1 to $\underline{\mathbf{C}}_{-t}$ to obtain a set of $I$-kept rectangles which we call $\underline{\mathbf{K}}(t,\mathfrak{b}) \subset \underline{\mathbf{C}}_{-t}$ such that the configuration $\underline{\zeta}_t$, defined by $\underline{\zeta}_t(i) := |\{(\mathsf{R},U) \in \underline{\mathbf{K}}(t,\mathfrak{b}) : \gamma(\mathsf{R})_0 = i\}|$, has law $T_t(\nu_\alpha)$. The proof follows the same lines as that of Lemma 3.4. Note that the $U$ variables played no role in building $\underline{\mathbf{K}}(t,\mathfrak{b})$ from $\underline{\mathbf{C}}_{-t}$.

We deal now with the conditioned process starting at time $-t$ with the Gibbs measure $\mu_G$. First, we need to introduce a new type of parent. We say that $(\mathsf{R}',U') \in \mathbf{C}_U$ is a $\mu$-parent of $(\mathsf{R},U)$ if

(3.23) $$U > \frac{\alpha_i}{\rho} \quad \text{with } i = \gamma(\mathsf{R})_{-t},$$

$$\mathrm{epoch}(\mathsf{R}') \ni \mathrm{birth}(\mathsf{R}) \text{ and } |\gamma(\mathsf{R}')_{-t} - i| < R.$$

Note that when $U < \alpha_i/\rho$, then the birth of the rectangle $(\mathsf{R},U)$ is certain since the minimum of $c_i(\eta)$ over $\eta$ is larger than $\alpha_i/\rho$ by hypotheses.

REMARK 3.7. Actually, $\mu$-parents depend only on the sequence $\{\alpha_i, i \in \mathbb{Z}^d\}$ introduced in Proposition 1.4. Thus, we could have used the notation $\alpha$-parents. We keep the name of the initial measure $\mu$ to distinguish with the previous case where the initial measure is $\nu_\rho$.

Let $\mathbf{B}_1^{\mathsf{R},U}(I,\mu)$ be the first generation of both $\mu$-parents and $I$-parents [as defined in Section 3.2.1, but considered as couples $(\mathsf{R},U)$]. Let $\mathbf{B}_n^{\mathsf{R},U}(I,\mu)$ and $\mathbf{B}^{\mathsf{R},U}(I,\mu)$ be respectively the $n$th generation and the clan of $(\mathsf{R},U)$ defined by

(3.24) $$\mathbf{B}_n^{\mathsf{R},U}(I,\mu) = \bigcup_{(\mathsf{R}',U') \in \mathbf{B}_{n-1}^{\mathsf{R},U}(I,\mu)} \mathbf{B}_1^{\mathsf{R}',U'}(I,\mu) \quad \text{and}$$

$$\mathbf{B}^{\mathsf{R},U}(I,\mu) = \bigcup_{n \geq 1} \mathbf{B}_n^{\mathsf{R},U}(I,\mu) \cup \{(\mathsf{R},U)\}.$$

LEMMA 3.8. *Fix $t > 0$ and let $I = [-t, 0]$. For $\rho < \rho_c(\beta_1)$ given in* (1.6), *and $\mu_G$ given in Proposition* 1.4, *we have that* $\mathbf{B}^{\mathsf{R},U}(I, \mu)$ *is $\mathbb{Q}$-almost surely finite.*

Lemma 3.8 is proved in Section 6.3.

*The $(I, \mu)$-trimming algorithm.* Fix a time $\mathfrak{b}$ and let $(\mathsf{R}, U) \in \mathbf{C}_U$ be alive at $\mathfrak{b}$.

- Let $(\mathsf{R}_1, U_1)$ be the eldest element of $\mathbf{B}^{\mathsf{R},U}(I, \mu)$. If $\Gamma_1 := \{\gamma(\mathsf{R}_1)\} \in \mathcal{A}_I$, or
$$U_1 > c_i(\eta_{-t}(\Gamma_1)) \qquad \text{with } i = \gamma(\mathsf{R}_1)_{-t},$$
then $(\mathsf{R}_1, U_1)$ is deleted, else it is kept.
- Assume we have labeled $C_n := \{(\mathsf{R}_1, U_1), \dots, (\mathsf{R}_n, U_n)\}$, the eldest $n$ elements of the clan. If
$$\Gamma_n := \{\gamma(\mathsf{R}) : (\mathsf{R}, U) \in C_n, \text{ is kept}, \text{epoch}(\mathsf{R}) \ni \text{birth}(\mathsf{R}_{n+1})\}$$
$$\cup \{\gamma(\mathsf{R}_{n+1})\} \in \mathcal{A}_I,$$
or
$$U_{n+1} > c_i(\eta_{-t}(\Gamma_n)) \qquad \text{with } i = \gamma(\mathsf{R}_{n+1})_{-t},$$
then $(\mathsf{R}_{n+1}, U_{n+1})$ is deleted, else it is kept.
- Stop the labeling once all elements in the clan are labeled.

Repeating this procedure with all elements of $\mathbf{C}_U$ alive at $\mathfrak{b}$, we build the set of kept rectangles denoted by $\mathbf{K}(I, \mu, \mathfrak{b})$, and we define $\Gamma^{\mathfrak{b}}_{I,\mu} := \{\gamma(\mathsf{R}) : (\mathsf{R}, U) \in \mathbf{K}(I, \mu, \mathfrak{b})\}$. We omit the easy proofs of the following proposition (and refer the interested reader to similar arguments in [7], proof of Theorem 3.1, and [8], Appendix A, proof of Theorem 1).

PROPOSITION 3.9. *Fix $t > 0$ and let $I = [-t, 0]$. The process $\{\Gamma^{\mathfrak{b}}_{I,\mu} : \mathfrak{b} \in \mathbb{R}\}$ is stationary and it is Markov with generator $\mathcal{L}^{I,\mu}_{\mathrm{bd}}$ defined by*

$$
\begin{aligned}
\mathcal{L}^{I,\mu}_{\mathrm{bd}} f(\Gamma) = &\sum_{i \in \mathbb{Z}^d} \int P_{-t,i}(d\gamma) c_i(\eta_{-t}(\Gamma)) \\
&\times \mathbb{1}\{\Gamma \cup \{\gamma\} \in \mathcal{A}^c_I\}[f(\Gamma \cup \{\gamma\}) - f(\Gamma)] \\
&+ \sum_{\gamma \in \Gamma} [f(\Gamma \setminus \{\gamma\}) - f(\Gamma)].
\end{aligned}
$$
(3.25)

*Furthermore, the law of $\Gamma^{\mathfrak{b}}_{I,\mu}$ is the measure induced in $\mathcal{G}$ by the AIRW process starting with $\mu_G$ at time $-t$ and conditioned on not hitting $\mathcal{A}$ in the interval $I$. The time-slice $\eta_0(\Gamma^{\mathfrak{b}}_{I,\mu})$ has law $T_t(\mu_G)$.*

**4. Estimating discrepancies.** Henceforth, we assume $d \geq 3$. Let $\beta_d$ be given in (A.9) and $\rho_c(\beta_1)$ be given in (1.6). We fix $\beta_1 < \beta_d$ and $\rho < \rho_c(\beta_1)$. Also, we fix an arbitrary $\mathfrak{b} \in \mathbb{R}$.

4.1. *Discrepancies between* $T_t(\nu_\rho)$ *and* $T_\infty(\nu_\rho)$. We call, for notational convenience, $\zeta_t := \eta_0(\Gamma^{\mathfrak{b}}_{[-t,0]})$ and $\zeta_\infty := \eta_0(\Gamma^{\mathfrak{b}}_{(-\infty,0]})$. This is a coupling of $T_t(\nu_\rho)$ and $T_\infty(\nu_\rho)$ built in Section 3.2.1 as a deterministic function of $\mathbf{C}$. Also, for a realization of $\mathbf{C}$, we define, for $i \in \mathbb{Z}^d$,

(4.1) $\quad \underline{\xi}(i) := |\{\mathsf{R} \in \mathbf{C} : \gamma(\mathsf{R})_0 = i, \Sigma(\mathsf{R}) \cap (-\infty, 0] = \varnothing, \text{epoch}(\mathsf{R}) \ni \mathfrak{b}\}|.$

The main result of this section is the following.

LEMMA 4.1. *The variable* $\underline{\xi} \prec \zeta_t \wedge \zeta_\infty$, *and there is an explicit configuration* $\bar{\xi}_t$ *satisfying:*

- $\bar{\xi}_t$ *and* $\underline{\xi}$ *are independent.*
- $|\zeta_t - \zeta_\infty| \prec \bar{\xi}_t.$

*Furthermore, there are positive constants $C_1, C_2$ such that, for any site $i \in \mathbb{Z}^d$,*

$$\begin{aligned}
E[\bar{\xi}_t(i)] &\leq \rho P^*_{0,i}(t < H_0 < \infty) + C_2 P^*_{0,i}(H_0 < \infty) e^{-\beta_1 t} \\
&\quad + C_1 \sum_j p(0,j) \int_0^t P_{0,j}(\gamma_{t-s} = i, H_0 = \infty) e^{-\beta_1 s} \, ds.
\end{aligned} \tag{4.2}$$

The following corollary combines Lemma 4.1 and Lemma A.4 of the Appendix.

COROLLARY 4.2. *There is a constant $C(\beta_1)$ such that*

(4.3) $$\sum_{i \in \mathbb{Z}^d} P_{0,i}(H_0 < \infty) E[\bar{\xi}_t(i)] \leq C(\beta_1) e^{-\beta_1 t}$$

*and*

(4.4) $$\lim_{t \to \infty} \sum_{i \in \mathbb{Z}^d} P^*_{0,i}(H_0 < \infty) E[\bar{\xi}_t(i)] = 0.$$

PROOF OF LEMMA 4.1. First, we characterize $\underline{\xi}$. For any $I \subset (-\infty, 0]$, the rectangles making up $\underline{\xi}$ have no $I$-parents, and are always $I$-kept so that $\underline{\xi} \prec \zeta_t \wedge \zeta_\infty$. Note that in (4.1), $\underline{\xi}$ is distributed as $\hat{T}_\infty(\nu_\rho)$, with the notation of Section 3.3 with $\Lambda = \{0\}$. Thus, its law is a product of Poisson laws with marginal at site $j$ of mean $\rho P^*_{0,j}(H_0 = \infty)$.

Second, we consider the discrepancies. Note that $\mathbf{A}^{\mathsf{R}}[-t, 0] \subset \mathbf{A}^{\mathsf{R}}(-\infty, 0]$, and that if $\mathbf{A}^{\mathsf{R}}(-\infty, 0] = \mathbf{A}^{\mathsf{R}}[-t, 0]$, then $\mathsf{R}$ would not be a discrepancy since it would have the same fate under the $I$-trimming algorithm for both $I = (-\infty, 0]$ and $I = [-t, 0]$. We now define

(4.5) $$\tilde{\mathbf{C}} := \{\mathsf{R} \in \mathbf{C} : \Sigma(\mathsf{R}) \cap (-\infty, 0] \neq \varnothing\},$$

and observe that

$$|\zeta_t - \zeta_\infty|(i) = \left| \sum_{\mathsf{R} \in \tilde{\mathbf{C}}} \mathbb{1}\{\gamma(\mathsf{R})_0 = i\} \right.$$

$$\left. \times (\mathbb{1}\{\mathsf{R} \in \mathbf{K}([-t, 0], \mathfrak{b})\} - \mathbb{1}\{\mathsf{R} \in \mathbf{K}((-\infty, 0], \mathfrak{b})\}) \right|$$

(4.6)
$$\leq \bar{\xi}_t(i) := \sum_{\mathsf{R} \in \tilde{\mathbf{C}}} \mathbb{1}\{\gamma(\mathsf{R})_0 = i, \operatorname{epoch}(\mathsf{R}) \ni \mathfrak{b}\}$$

$$\times \mathbb{1}\{\mathbf{A}^\mathsf{R}[-t, 0] \neq \mathbf{A}^\mathsf{R}(-\infty, 0]\}.$$

Also, since $\bar{\xi}_t$ is a function of $\tilde{\mathbf{C}}$, and $\underline{\xi}$ is a function of $\tilde{\mathbf{C}}^c$, $\bar{\xi}_t$ and $\underline{\xi}$ are independent. □

Finally, we divide the rectangles contributing to $\bar{\xi}_t$ into three collections that we treat separately. We need more notation: let $\bar{s}(\mathsf{R}) = \sup\{\boldsymbol{\sigma}(\mathsf{R})\}$, and $\underline{s}(\mathsf{R}) = \inf\{\boldsymbol{\sigma}(\mathsf{R})\}$, and define the following:

- $\mathbf{C}^0 = \{\mathsf{R} \in \tilde{\mathbf{C}} : \boldsymbol{\sigma}(\mathsf{R}) \subset (-\infty, -t]\}$, with which we associate $\xi_t^0$. That is to say

$$\forall i \in \mathbb{Z}^d \quad \xi_t^0(i) = |\{\mathsf{R} \in \mathbf{C}^0 : \gamma(\mathsf{R})_0 = i, \operatorname{epoch}(\mathsf{R}) \ni \mathfrak{b}\}|.$$

- $\mathbf{C}^1 = \{\mathsf{R} \in \tilde{\mathbf{C}} : \bar{s}(\mathsf{R}) \in ]-t, 0], \mathbf{A}^\mathsf{R}[-t, 0] \neq \mathbf{A}^\mathsf{R}(-\infty, 0]\}$. We associate $\xi_t^1$ to $\mathbf{C}^1$.
- $\mathbf{C}^2 = \{\mathsf{R} \in \tilde{\mathbf{C}} : \bar{s}(\mathsf{R}) > 0, \underline{s}(\mathsf{R}) < 0, \mathbf{A}^\mathsf{R}[-t, 0] \neq \mathbf{A}^\mathsf{R}(-\infty, 0]\}$. We associate $\xi_t^2$ to $\mathbf{C}^2$.

*Case* $\xi_t^0$. The rectangles which participate to $\mathbf{C}^0$ are associated with trajectories drawn independently at $t = 0$, which do not meet site $\{0\}$ during $[-t, 0]$. Since $\mathbb{Q}$ consists in drawing at each site $i \in \mathbb{Z}^d$ a Poisson process of intensity $\rho$ whose time-realizations are marked with a trajectory drawn from $P_{0,i}$, it is obvious that

$$E[\xi_t^0(i)] = \rho P_{0,i}^*(t < H_0 < \infty).$$

*Case* $\xi_t^1$. For $\mathsf{R} \in \mathbf{C}^1$, since $\mathbf{A}^\mathsf{R}[-t, 0] \neq \mathbf{A}^\mathsf{R}(-\infty, 0]$, we have $W(\mathsf{R}) \cap (-\infty, -t] \neq \varnothing$, by Remark 3.5. Thus, estimates on the width of a clan, obtained in Section 6.1, will play a key role here. But first, it is convenient to give an alternative construction of the marked-Poisson process corresponding to $\mathbf{C}^1$, by marking the last visit time to 0, that is, by considering

$$(4.7) \qquad \{(\bar{s}(\mathsf{R}), \mathrm{birth}(\mathsf{R}), \mathrm{death}(\mathsf{R})), \mathsf{R} \in \mathbf{C}^1\},$$

and by partitioning this set in terms of first exit sites from $\{0\}$: that is, the set of $j \neq 0$ such that at $t = \bar{s}(\mathsf{R})$, $\gamma(\mathsf{R})_{t-} = 0$ and $\gamma(\mathsf{R})_{t+} = j$. As $j$ runs over $\{j : p(0, j) > 0\}$, we obtain independent point processes with respective intensity $g^j(\bar{s}) \, d\bar{s} \, db \, e^{-l} \, dl$ with

$$(4.8) \qquad g^j(\bar{s}) := \mathbb{1}\{-t < \bar{s} < 0\} \rho p(0, j) P_{0,j}(H_0 = \infty).$$

*Marking the point process.* For each $j$ of $\{j : p(0, j) > 0\}$, we mark each point $(\bar{s}, b, l)$ by a trajectory $\tilde{\gamma}$ made up by concatenating two trajectories in the following way:

- $\gamma^+$ is drawn from $\mathbb{1}\{H_0 = \infty\} \, dP_{0,j}(\gamma^+)/P_{0,j}(\{H_0 = \infty\})$ (this is what happens after the last visit to 0 when the exit is from $j$).
- $\gamma^-$ is drawn from $dP^*_{0,0}$ (but we take trajectories which are left continuous with a right limit, so that the time-reversed trajectory has the correct shape).

Now, for $s < \bar{s}$, $\tilde{\gamma}_s = \gamma^-_{\bar{s}-s}$, whereas for $s \geq \bar{s}$, $\tilde{\gamma}_s = \gamma^+_{s-\bar{s}}$. Thus, $\tilde{\gamma}$ is drawn from $\sum_i dP_{0,i}$ conditioned on making the jump $0 \to j$ at time $\bar{s}$, and never visiting 0 after time $\bar{s}$. Note also that $\boldsymbol{\sigma}(\tilde{\gamma}) = \boldsymbol{\sigma}(\gamma^-)$ is independent of $\gamma^+$. For each $j$ of $\{j : p(0, j) > 0\}$, we denote the above mentioned marked-point process by

$$(4.9) \qquad \{N^j(\bar{s}, \tilde{\gamma}, b, l), \bar{s} \in \mathbb{R}, \tilde{\gamma} \in \mathcal{D}(\mathbb{R}, \mathbb{Z}^d), b \in \mathbb{R}, l \in \mathbb{R}^+\}.$$

It is clear that the corresponding rectangle process

$$\left\{ \tilde{\gamma} \times [b, b+l[ : (\bar{s}, \tilde{\gamma}, b, l) \in \mathrm{support}\left(\bigcup_j N^j\right) \right\}$$

has the same law as $\mathbf{C}^1$.

The point of such a representation is that, conditioned on $\bar{s}$, we have independence of $\{\gamma^+_{|\bar{s}|} = i\}$ from the width of $\mathsf{R} = (\gamma^-, \gamma^+) \times [b, b+l[$, $W(\mathsf{R})$, defined in (3.15) to be the union of $\sigma(\gamma(\mathsf{R}')) := \bar{s}(\gamma(\mathsf{R}')) - \underline{s}(\gamma(\mathsf{R}'))$, where $\mathsf{R}'$ runs over $\mathbf{A}^\mathsf{R}((-\infty, 0])$. Now, since $W(\mathsf{R})$ depends only on $\{(\bar{s}(\mathsf{R}'), \sigma(\mathsf{R}')) : \mathsf{R}' \in \mathbf{A}^\mathsf{R}((-\infty, 0])\}$, we can further simplify the description of the processes $\{N^j\}$. Thus, we consider the *projected marked processes* made up of

$$(4.10) \qquad \{x = (\bar{s}, \sigma(\tilde{\gamma}), b, l) : (\bar{s}, \tilde{\gamma}, b, l) \in \mathrm{support}(N^j)\},$$

that we still call $N^j$ for convenience, and we denote by $N_p := \bigcup \{N^j : j$ such that $p(0,j) > 0\}$. For $x \in \text{support}(N_p)$, we denote the width of the corresponding rectangle by $W(x)$. Now, we have

$$\xi_t^1(i) \leq \sum_{j\,:\,p(0,j)>0} \int \mathbb{1}\{\gamma_{|\bar{s}|} = i\} \mathbb{1}\{b \leq \mathfrak{b} < b + l\}$$

(4.11)
$$\times \mathbb{1}\{W(\bar{s}, \sigma, b, l) \ni \bar{s} + t\} \, dN^j(\bar{s}, \sigma, b, l)$$

$$\times \frac{\mathbb{1}\{H_0(\gamma) = \infty\} \, dP_{0,j}(\gamma)}{P_{0,j}(H_0 = \infty)}.$$

By taking expectation and performing an obvious change of variables,

$$E[\xi_t^1(i)] \leq \rho \sum_j p(0,j) \int_{\mathbb{R}} db \int_0^\infty dl \, e^{-l} \int_0^t ds \int g_\sigma(\sigma) \, d\sigma$$

(4.12)
$$\times \mathbb{1}\{b \leq \mathfrak{b} < b + l\}$$

$$\times P_{0,j}(\gamma_s = i, H_0 = \infty) P(W(0, \sigma, b, l) \ni t - s).$$

We show in Section 6.1 that $P(W(0, \sigma, b, l) \ni t - s)$ is independent of $b$ and $l$, and a simple computation yields

$$\forall \mathfrak{b} \in \mathbb{R} \qquad \int_{\mathbb{R}} db \int_0^\infty dl \, e^{-l} \mathbb{1}\{b \leq \mathfrak{b} < b + l\} = 1.$$

Now, by Lemma 6.1, there is a positive (explicit) constant $C_1$ such that

(4.13) $$E[\xi_t^1(i)] \leq C_1 \sum_j p(0,j) \int_0^t P_{0,j}(\gamma_s = i, H_0 = \infty) e^{-\beta_1(t-s)} \, ds.$$

*Case $\xi_t^2$.* We need here to condition on both $\bar{s}$ and $\sigma$ to obtain independence of the width of the clan, and of $\{\tilde{\gamma}_0(\mathsf{R}) = i\}$. Thus, we consider the point process

$$\{(\bar{s}(\mathsf{R}), \underline{s}(\mathsf{R}), \text{birth}(\mathsf{R}), \text{death}(\mathsf{R})), \mathsf{R} \in \tilde{\mathbf{C}}\},$$

which we partition into the last site before hitting $\{0\}$, say, $j \neq 0$, and the first exit site from $\{0\}$, say, $j' \neq 0$. Proceeding similarly as for $\xi_t^1$, for each $j$ and $j'$, with $p(j,0) > 0$ and $p(0,j') > 0$, we consider the point processes on $\{(\bar{s}, \underline{s}, b, l) \in \mathbb{R} \times \mathbb{R} \times \mathbb{R} \times \mathbb{R}^+\}$ with mean measure

(4.14)
$$\mathbb{1}\{\bar{s} > 0 > \underline{s}\} g_\sigma(\bar{s} - \underline{s}) \rho p(j, 0) P_{0,j}^*(H_0 = \infty) p(0, j')$$

$$\times P_{0,j'}(H_0 = \infty) \, d\bar{s} \, d\underline{s} \, db \, e^{-l} \, dl,$$

where $g_\sigma$ is the density of the law of $\bar{s} - \underline{s}$; in Lemma A.3 of the Appendix, we bound $g_\sigma$.

We denote by $\{N_j^{j'}, j, j'\}$ the marked-point process obtained by attaching to the previous points trajectories denoted by $\tilde{\gamma}$ obtained as follows. Fix $\bar{s} > 0$ and $\underline{s} < 0$, and $\sigma = \bar{s} - \underline{s}$:

- Draw $\gamma^m$ and $\gamma^-$ respectively from

$$\frac{\mathbb{1}\{\gamma^m_\sigma = 0\} \, dP_{0,0}(\gamma^m)}{P_{0,0}(\gamma_\sigma = 0)} \quad \text{and} \quad \frac{\mathbb{1}\{H_0(\gamma^-) = \infty\} \, dP^*_{0,j}(\gamma^-)}{P^*_{0,j}(\{H_0 = \infty\})}.$$

- Draw $\gamma^+$ as in the previous case.
- Concatenate the time-reversed $\gamma^-$ before time $\underline{s}$, $\gamma^m$ between time $\underline{s}$ and $\bar{s}$, and $\gamma^+$ after time $\bar{s}$, to obtain $\tilde{\gamma}$.

Then, we have the bound

$$\xi_t^2(i) \le \sum_j \sum_{j'} \int \mathbb{1}\{\gamma^m_{|\underline{s}|} = i\} \mathbb{1}\{b \le \mathfrak{b} < b+l\}$$
(4.15)
$$\times \mathbb{1}\{W(\bar{s},\sigma,b,l) \ni t + \bar{s}\} \, dN^{j'}_j(\bar{s}, \underline{s}, \tilde{\gamma}, b, l).$$

Fixing $\bar{s}$ and $\underline{s}$, we have the conditional independence of $\{\gamma^m_{|\underline{s}|} = i\}$ and $\{W(\bar{s},\sigma,b,l) \ni t + \bar{s}\}$. Thus, after time-shifting the variables $\underline{s}$ and $\bar{s}$, integrating over $b$ and $l$, and summing the $j$ and $j'$, we obtain

(4.16)
$$\begin{aligned}
E[\xi_t^2(i)] &\le \rho \int_t^\infty d\bar{s} \int_{-\infty}^t d\underline{s} \, P_{0,0}(\gamma_{t-\underline{s}} = i) g_\sigma(\bar{s} - \underline{s}) \\
&\quad \times P(W(0, \bar{s} - \underline{s}, b, l) \ni \bar{s}) \\
&\le \rho \int_{-\infty}^t d\underline{s} \, P^*_{0,i}(H_0 < t - \underline{s}) \int_t^\infty d\bar{s} \, g_\sigma(\bar{s} - \underline{s}) \\
&\quad \times P(W(0, \bar{s} - \underline{s}, b, l) \ni \bar{s}).
\end{aligned}$$

We use that $t \mapsto P(W(0,\sigma,b,l) \ni t)$ is decreasing in $t \ge 0$ to obtain

(4.17)
$$E[\xi_t^2(i)] \le \rho P^*_{0,i}(H_0 < \infty) \Bigg( \int_0^t ds \int_0^\infty g_\sigma(\sigma) P(W(0,\sigma,b,l) \ni t) \, d\sigma + \int_t^\infty ds \int_s^\infty g_\sigma(\sigma) \, d\sigma \Bigg).$$

Finally, using Lemmas 6.1 and A.3, there is a constant $C_2$ such that

(4.18)
$$E[\xi_t^2(i)] \le C_2 P^*_{0,i}(H_0 < \infty) e^{-\beta_1 t}.$$

4.2. *Discrepancies between $T_t(\nu_\rho)$ and $T_t(\mu_G)$.* We fix $\mathfrak{b}$ and $t > 0$ and consider henceforth $I = [-t, 0]$. We call, for notational convenience, $\zeta_t := \eta_0(\Gamma^\mathfrak{b}_I)$ and $\zeta_t^\mu := \eta_0(\Gamma^\mathfrak{b}_{I,\mu})$. We recall that we have introduced in Section 3.4 the set of *rectangles* $\mathbf{C}_U$, whose elements are of the form $(\mathsf{R}, U)$ with $\mathsf{R} \in \mathbf{C}$ and $U$ are uniform variables in $[0, 1]$.

We first build a configuration $\underline{\xi}_t$ coming from rectangles of $K(I, \mathfrak{b}) \cap K(I, \mu, \mathfrak{b})$. We say that a rectangle $(\mathsf{R}, U) \in \mathbf{C}_U$ is *good* if it is alive at $\mathfrak{b}$ and

(i) $U \leq \alpha_i/\rho$ where $i = \gamma(\mathsf{R})_{-t}$, and (ii) $\Sigma(\mathsf{R}) \cap I = \varnothing$. Thus, a *good* rectangle has neither $\mu$-parents nor $I$-parents. We define, for each $i \in \mathbb{Z}^d$,

(4.19)
$$\underline{\xi}_t(i) := |\{(\mathsf{R}, U) \ good : \gamma_0(\mathsf{R}) = i\}|.$$

LEMMA 4.3. *For any $t > 0$, the configuration $\underline{\xi}_t \prec \zeta_t \wedge \zeta_t^\mu$, and there is an explicit $\bar{\xi}_t$ satisfying the following:*

- *$\bar{\xi}_t$ and $\underline{\xi}_t$ are independent.*
- *$|\zeta_t - \zeta_t^\mu| \prec \bar{\xi}_t$.*
- *$\underline{\xi}_t$ is distributed according to a product of Poisson law, that we denote by $\nu_{\beta^t}$. Moreover, $\{\nu_{\beta^t}, t \geq 0\}$ have densities with respect to $\nu_\rho$ which are uniformly bounded in $L^2(\nu_\rho)$.*

*Furthermore, there is a constant $C$ such that, for any site $i \in \mathbb{Z}^d$,*

(4.20)
$$E[\bar{\xi}_t(i)] \leq C \sum_{j \neq 0} \int \int_0^t \rho p(0,j) P_{0,j}(\gamma_{t-s} = i, H_0 = \infty) e^{-\beta_1 s} \, ds$$
$$+ C P_{0,i}^*(H_0 < \infty) e^{-\beta_1 t}.$$

As a corollary of Lemmas 4.3 and A.4, we have the following.

COROLLARY 4.4. *There is a positive constant $C$ such that*

(4.21)
$$\sum_{i \in \mathbb{Z}^d} P_{0,i}(H_0 < \infty) E[\bar{\xi}_t(i)] \leq C e^{-\beta_1 t}$$

*and*

(4.22)
$$\lim_{t \to \infty} \sum_{i \in \mathbb{Z}^d} P_{0,i}^*(H_0 < \infty) E[\bar{\xi}_t(i)] = 0.$$

PROOF OF LEMMA 4.3. Since the proof is close to the proof of Lemma 4.1, we mainly focus on the differences.

The first step is to characterize $\underline{\xi}_t$. Note that the *good* rectangles have no $(I, \mu)$-parents and, thus, $\underline{\xi}_t \prec \zeta_t \wedge \zeta_t^\mu$. Now, we consider the law of $\underline{\xi}_t$: at time $-t$, the trajectories are picked up with rate $\alpha_i$ at site $i$, and only those with $H_0 \circ \theta_t > t$ are kept. Thus, $\{\underline{\xi}_t(j), i \in \mathbb{Z}^d\}$ are independent Poisson variables of intensity

$$\beta_j^t = \sum_{i \in \mathbb{Z}^d} \alpha_i P_{0,i}(\gamma_t = j, H_0 > t) \qquad \text{at site } j \in \mathbb{Z}^d.$$

Now, by reversing time and using (1.9),

$$1 - \frac{\beta_j^t}{\rho} = 1 - \sum_{i \in \mathbb{Z}^d} \frac{\alpha_i}{\rho} P_{0,j}^*(\gamma_t = i, H_0 > t)$$

$$= 1 - P_{0,j}^*(H_0 > t) + \sum_{i \in \mathbb{Z}^d} \left(1 - \frac{\alpha_i}{\rho}\right) P_{0,j}^*(\gamma_t = i, H_0 > t)$$

(4.23)

$$\leq P_{0,j}^*(H_0 < \infty) + C_\alpha \sum_{i \in \mathbb{Z}^d} P_{0,j}^*(\gamma_t = i, H_0 > t) P_{0,i}^*(H_0 < \infty)$$

$$\leq (1 + C_\alpha) P_{0,j}^*(H_0 < \infty).$$

Thus, the densities of the laws of $\{\underline{\xi}_t, t \geq 0\}$ are uniformly bounded in $L^2(\nu_\rho)$ as soon as $d \geq 3$, by Remark 1.5.

The second step is to deal with discrepancies. A rectangle $(\mathsf{R}, U)$ can make up a discrepancy between $\zeta_t$ and $\zeta_t^\mu$ if either it has $\mu$-parents [and, therefore, $U > \alpha_i/\rho$, where $i = \gamma(\mathsf{R})_{-t}$], or one of his $I$-parents has a $\mu$-parent. Since we only need to overcount the discrepancies, we introduce the following subset $\mathbf{C}_{\text{bad}}$ of $\mathbf{C}_U$ of *bad* rectangles. A rectangle $(\mathsf{R}, U)$ is *bad* when $U > \alpha_i/\rho$, where $i = \gamma(\mathsf{R})_{-t}$. Then, a rectangle is susceptible of being a discrepancy if one of its $I$-ancestors is *bad*. Thus,

(4.24) $\overline{\xi}_t(i) := |\{(\mathsf{R}, U) \in \mathbf{C}_U : \mathbf{A}^\mathsf{R}(I) \cap \mathbf{C}_{\text{bad}} \neq \varnothing, \gamma(\mathsf{R})_0 = i, \text{epoch}(\mathsf{R}) \ni \mathfrak{b}\}|.$

Thus, it is clear that $\underline{\xi}_t$ and $\bar{\xi}_t$ are independent.

Following Section 4.1, we divide the rectangles contributing to $\bar{\xi}_t$ into two sets:

- $\mathbf{C}^1 = \{(\mathsf{R}, U) \in \mathbf{C}_U : \mathbf{A}^\mathsf{R}(I) \cap \mathbf{C}_{\text{bad}} \neq \varnothing, \bar{s}(\mathsf{R}) \in [-t, 0]\}$, with which we associate $\xi_t^1$.
- $\mathbf{C}^2 = \{(\mathsf{R}, U) \in \mathbf{C}_U : \mathbf{A}^\mathsf{R}(I) \cap \mathbf{C}_{\text{bad}} \neq \varnothing, \bar{s}(\mathsf{R}) > 0 > \underline{s}(\mathsf{R})\}$, with which we associate $\xi_t^2$.

*Case $\xi_t^1$.* As in the preceding section, we consider the projected marked-Poisson process $N_p$ given in (4.10) which will correspond to the $I$-ancestors. As in the proof of Lemma 4.1, we give an alternative construction of $\mathbf{C}_{\text{bad}}$. Thus, we consider

$$\{(\bar{s}(\mathsf{R}), \text{birth}(\mathsf{R}), \text{death}(\mathsf{R})), (\mathsf{R}, U) \in \mathbf{C}_{\text{bad}}\},$$

and we partition this set in terms of exit sites $j \neq 0$ with $\gamma(\mathsf{R})_{t-} = 0$ and $\gamma(\mathsf{R})_{t+} = j$ with $t = \bar{s}(\mathsf{R})$. As $j$ covers $\{j : p(0, j) > 0\}$, we obtain independent point processes with respective intensity $\tilde{g}^j(\bar{s}) \, d\bar{s} \, db \, e^{-l} \, dl$. We first show that, for any $\beta_1 < \beta_d$, there is a constant $c_0(\beta_1)$ such that

(4.25) $$\sum_j \tilde{g}^j(\bar{s}) \leq \rho c_0(\beta_1) \mathbb{1}\{-t < \bar{s} < 0\} e^{-\beta_1(\bar{s}+t)},$$

where $c_0(\beta_1) := c_1(\beta_1)C_\alpha < \infty$ and $M(\beta_1)$ is given in Lemma A.4(v). First, when $\bar{s} + t \geq 0$, it is easy to see that

$$\tilde{g}^j(\bar{s}) = p(0,j)P_{0,j}(H_0 = \infty) \sum_{i \in \mathbb{Z}^d} (\rho - \alpha_i)P^*_{0,0}(\gamma_{\bar{s}+t} = i)$$

(4.26)
$$\leq \rho C_\alpha p(0,j)P_{0,j}(H_0 = \infty) \sum_{i \in \mathbb{Z}^d} P^*_{0,0}(\gamma_s = i)P^*_{0,i}(H_0 < \infty),$$

where we used $1 - \alpha_i/\rho \leq C_\alpha P^*_{0,i}(H_0 < \infty)$. Now, (4.25) follows from the estimate (v) of Lemma A.4. Now, we mark each point $(\bar{s}, b, l)$, as in the paragraph following (4.8), by the time spent at site 0. Thus, since we need to overcount discrepancies, we overcount the number of *bad* points by introducing, independently of $N_p$, a Poisson process $\tilde{N}$ on $E^+ := \{x = (s, \sigma, b, l) \in \mathbb{R}^+ \times \mathbb{R}^+ \times \mathbb{R} \times \mathbb{R}^+\}$, whose intensity measure has density

(4.27) $$\rho c_0(\beta_1) \exp(-\beta_1 s) \, ds \, g_\sigma(\sigma) \, d\sigma \, db \, \exp(-l) \, dl.$$

Note that the variable $s$ corresponds to $\bar{s} + t \geq 0$. We denote by $\tilde{m}$ the support of $\tilde{N}$. We call the points of $\tilde{m}$ the *bad* points. For a point $x = (s, \sigma, b, l)$ of $N_p$, we introduce $K(x)$, which is the event that $x$ has an $I$-ancestor which is *bad*. In Section 6.2 we estimate the probability of $K(x)$.

In this section $K(x)$ plays the role that the width of the clan played in Section 4.1. Here also, conditioned on $\bar{s}$, $K(x)$ is independent of $\gamma^+_{|\bar{s}|} = i$. Thus, using the point processes $N^j$ defined in (4.10), we obtain

$$\xi^1_t(i) \leq \sum_{j:p(0,j)>0} \int \mathbb{1}\{\gamma_{|\bar{s}|} = i\}\mathbb{1}\{b \leq \mathfrak{b} \leq b + l\}$$

(4.28)
$$\times \mathbb{1}\{K(\bar{s}, \sigma, b, l)\} \, dN^j(\bar{s}, \sigma, b, l)$$
$$\times \frac{\mathbb{1}\{H_0(\gamma) = \infty\} \, dP_{0,j}(\gamma)}{P_{0,j}(H_0 = \infty)}$$

and, after integration [see (4.12) for some intermediary steps],

$$E[\xi^1_t(i)] \leq \rho \sum_j p(0,j) \int_0^\infty g_\sigma(\sigma) \, d\sigma \int_0^t ds$$

(4.29)
$$\times P_{0,j}(\gamma_s = i, H_0 = \infty)P(K(t-s, \sigma, b, l)).$$

Finally, Lemma 6.2 is used to obtain inequality (4.13).

*Case $\xi^2_t$.* This case is also analogous to that of the previous section. Let $\{N^{j'}_j\}$ be independent Poisson processes (corresponding to $I$-ancestors) with intensities given in (4.14). As for (4.15), we have

$$\xi^2_t(i) \leq \sum_j \sum_{j'} \int \mathbb{1}\{\gamma^m_{|\underline{s}|} = i\}\mathbb{1}\{K(\bar{s}, \sigma, b, l)\}$$

(4.30)
$$\times \mathbb{1}\{b \leq \mathfrak{b} \leq b + l\} \, dN^{j'}_j(\bar{s}, \underline{s}, \tilde{\gamma}, b, l),$$

and, after integration,

$$(4.31) \quad E[\xi_t^2(i)] \leq \rho \int_t^\infty d\bar{s} \int_{-\infty}^t d\underline{s}\, g_\sigma(\bar{s} - \underline{s}) P_{0,0}(\gamma_{t-\underline{s}} = i) P(K(\bar{s}, \bar{s} - \underline{s}, b, l)).$$

Lemma 6.2 is used to obtain inequality (4.18) using similar arguments as in Section 4.1. □

## 5. Hitting times.

5.1. *Proof of Proposition* 1.1. We recall that in [3], principal Dirichlet eigenfunctions denoted by $u$ and $u^*$, respectively for $\mathcal{L}_W$ and $\mathcal{L}_W^*$, were shown to exist in $L^p(\nu_\rho)$ for any integer $p$, by considering limits of linear combinations of respectively

$$(5.1) \qquad u_t(\eta) := \frac{P_\eta(\tau > t)}{P_{\nu_\rho}(\tau > t)} \quad \text{and} \quad u_t^*(\eta) := \frac{P_\eta^*(\tau > t)}{P_{\nu_\rho}^*(\tau > t)}.$$

Note that $u_t = dT_t^*(\nu_\rho)/d\nu_\rho$ [and, similarly, $u_t^* = dT_t(\nu_\rho)/d\nu_\rho$]. Indeed, for any $\varphi \in L^2(\nu_\rho)$, we have

$$\int \varphi\, dT_t^*(\nu_\rho) = \frac{\int \bar{S}_t^*(\varphi) \mathbb{1}\{\mathcal{A}^c\}\, d\nu_\rho}{P_{\nu_\rho}(\tau > t)} = \frac{\int \varphi \bar{S}_t(\mathbb{1}\{\mathcal{A}^c\})\, d\nu_\rho}{P_{\nu_\rho}(\tau > t)} = \int \varphi u_t\, d\nu_\rho.$$

We recall also that in an $L^2(\nu_\rho)$-sense, $u$ and $u^*$ satisfy

$$\bar{S}_t(u) = e^{-\lambda t} u \quad \text{and} \quad \bar{S}_t^*(u^*) = e^{-\lambda t} u^*.$$

The key result in this section is the following.

LEMMA 5.1.   *Assume* $d \geq 3$. *For any* $\beta_1 < \beta_d$ *(where* $\beta_d$ *is given in Lemma* A.1*) and* $\rho < \rho_c(\beta_1)$, *there is an explicit number* $M(\beta_1)$ *such that, for any* $t \geq 0$,

$$(5.2) \qquad \sup_{s,t' \geq t} \left| \int u_s(u_t^* - u_{t'}^*)\, d\nu_\rho \right| \leq M(\beta_1) \exp(-\beta_1 t),$$

*and*

$$(5.3) \qquad \lim_{t \to \infty} \sup_{s,t' \geq t} \left| \int u_s(u_t - u_{t'})\, d\nu_\rho \right| = 0.$$

*The same result holds when* $\{u_t^*, t > 0\}$ *replaces* $\{u_t, t > 0\}$.

REMARK 5.2.   Propositions 1.2 and 1.1 are simple corollaries of Lemma 5.1. Indeed, first note that (5.3) implies that $\{u_t, t > 0\}$ and $\{u_t^*, t > 0\}$ are

$L^2(\nu_\rho)$-Cauchy sequences, since

$$\int (u_t - u_{t'})^2 \, d\nu_\rho = \int (u_t - u_{t'}) u_t \, d\nu_\rho - \int (u_t - u_{t'}) u_{t'} \, d\nu_\rho$$

$$\leq 2 \sup_{s,t' \geq t} \left| \int u_s (u_t - u_{t'}) \, d\nu_\rho \right| \stackrel{t \to \infty}{\longrightarrow} 0.$$

Thus, $u$ (resp. $u^*$) is the $L^2(\nu_\rho)$-limit of $\{u_t, t > 0\}$ (resp. of $\{u_t^*, t > 0\}$). Now, we take in (5.2) $s$ to infinity, then $t'$ to infinity to obtain

$$\left| \int u(u_t^* - u^*) \, d\nu_\rho \right| \leq M(\beta_1) e^{-\beta_1 t}.$$

By duality, $\int u u_t^* \, d\nu_\rho = \exp(-\lambda t)/P_{\nu_\rho}(\tau > t)$, and (1.7) follows with a constant $M(\beta_1) e^{\lambda t} P_{\nu_\rho}(\tau > t)/(\int u u^* \, d\nu_\rho)$. Now, since our event is decreasing $e^{\lambda t} P_{\nu_\rho}(\tau > t) \leq 1$ (see, e.g., [3]), and since $u$ and $u^*$ are decreasing, by FKGs inequality, $\int u u^* \, d\nu_\rho \geq \int u \, d\nu_\rho \int u^* \, d\nu_\rho = 1$. This yields (1.7).

PROOF OF LEMMA 5.1. First, note that

(5.4)
$$\left| \int u_s (u_t^* - u_{t'}^*) \, d\nu_\rho \right| \leq \left| \int u_s \, dT_t(\nu_\rho) - \int u_s \, dT_\infty(\nu_\rho) \right|$$
$$+ \left| \int u_s \, dT_{t'}(\nu_\rho) - \int u_s \, dT_\infty(\nu_\rho) \right|.$$

Thus, it is enough to treat

$$\int u_s \, dT_t(\nu_\rho) - \int u_s \, dT_\infty(\nu_\rho) = E[u_s(\zeta_t) - u_s(\zeta_\infty)],$$

where $(\zeta_t, \zeta_\infty)$ is the coupling of $dT_t(\nu_\rho)$ and $dT_\infty(\nu_\rho)$ introduced in Section 4.1. We recall that Lemma 4.1 establishes that there are two independent variables $\underline{\xi} \prec \zeta_t \wedge \zeta_\infty$ and $\bar{\xi}_t \succ |\zeta_t - \zeta_\infty|$ with $\underline{\xi} \sim \hat{\mu}_\rho$ [with the notation of (1.18) with $\Lambda = \{0\}$].

Now, we recall a simple observation. Since the process is monotone and $\mathcal{A}$ is increasing, $\eta \mapsto u_s(\eta)$ is decreasing for any $s$ and, by coupling the $\eta$-particles,

$$u_s(\eta) - u_s(\mathsf{A}_i^+ \eta) \leq u_s(\eta) P_{0,i}(H_0 < \infty).$$

By induction, this implies that

(5.5) $$|u_s(\zeta_t) - u_s(\zeta_\infty)| \leq u_s(\underline{\xi}) \sum_{i \in \mathbb{Z}^d} P_{0,i}(H_0 < \infty) \bar{\xi}_t(i).$$

A bound like (5.5) holds for $u_s^*$ with $P_{0,i}^*(H_0 < \infty)$ replacing $P_{0,i}(H_0 < \infty)$.

By taking the expectation of (5.5), and using independence of $\bar{\xi}_t$ and $\underline{\xi}$, we obtain

$$
\begin{aligned}
|E[u_s(\zeta_t) - u_s(\zeta_\infty)]| &\leq E[u_s(\underline{\xi})] \sum_{i \in \mathbb{Z}^d} E[\bar{\xi}_t(i)] P_{0,i}(H_0 < \infty) \\
&\leq C_1 \sum_{i \in \mathbb{Z}^d} E[\bar{\xi}_t(i)] P_{0,i}(H_0 < \infty),
\end{aligned}
\tag{5.6}
$$

with, by Cauchy–Schwarz and Proposition 1.7 [see (1.17) with $\Lambda = \{0\}$],

$$
\begin{aligned}
\sup_s E[u_s(\underline{\xi})] &= \sup_s \int u_s \, d\hat{\mu}_\rho \leq \sup_s(\|u_s\|_{\nu_\rho}) \left\| \frac{d\hat{\mu}_\rho}{d\nu_\rho} \right\|_{\nu_\rho} \\
&\leq C_1 := \left\| \frac{d\hat{\mu}_\rho^*}{d\nu_\rho} \right\|_{\nu_\rho} \left\| \frac{d\hat{\mu}_\rho}{d\nu_\rho} \right\|_{\nu_\rho} < \infty.
\end{aligned}
\tag{5.7}
$$

Finally, by (4.3) of Corollary 4.2, we deduce (5.2) from (5.6) with $M(\beta_1) := C_1 c_0(\beta_1)$.

Now, by a similar argument, we would have, for $s \geq t$,

$$
\begin{aligned}
\left| \int u_s^*(u_t^* \, d\nu_\rho - dT_\infty(\nu_\rho)) \right| &= |E[u_s^*(\zeta_t) - u_s^*(\zeta_\infty)]| \\
&\leq \left\| \frac{d\hat{\mu}_\rho}{d\nu_\rho} \right\|_{\nu_\rho}^2 \sum_{i \in \mathbb{Z}^d} E[\bar{\xi}_t(i)] P_{0,i}^*(H_0 < \infty).
\end{aligned}
\tag{5.8}
$$
$\square$

5.2. *Proof of Proposition* 1.4. We first express $f_t^* := dT_t(\mu_G)/d\nu_\rho$ in terms of the killed semi-group. For any function $\varphi$ in $L^2(\nu_\rho)$ with $\varphi|_{\mathcal{A}} \equiv 0$,

$$
(5.9) \quad \int \varphi \, dT_t(\mu_G) = \frac{\int \bar{S}_t(\varphi) \, d\mu_G}{P_{\mu_G}(\tau > t)} = \frac{\int \varphi \bar{S}_t^*(f_G) \, d\nu_\rho}{P_{\mu_G}(\tau > t)} \Longrightarrow f_t^* = \frac{\bar{S}_t^*(f_G)}{P_{\mu_G}(\tau > t)}.
$$

STEP 1. We prove that there is a sequence $\{\varepsilon_i, i \in \mathbb{Z}^d\}$ such that, for any $\eta$, and $i \in \mathbb{Z}^d$,

$$
(5.10) \qquad 0 \leq f_t^*(\eta) - f_t^*(\mathsf{A}_i^+ \eta) \leq \varepsilon_i f_t^*(\eta),
$$

with $\sum_i \varepsilon_i^2 < \infty$. It is only necessary to prove (5.10) for $\eta \mapsto \bar{S}_t^*(f_G)(\eta)$, since this inequality is homogeneous. Now, $\eta \mapsto \bar{S}_t^*(f_G)(\eta) = E_\eta^*[f_G(\eta_t) \mathbb{1}\{\tau > t\}]$ vanishes on $\mathcal{A}$ and is decreasing. Indeed, since the process is monotone, if $\eta \prec \zeta$, there is a coupling of the trajectories $(\eta., \zeta.)$ such that, almost surely, for all $t \geq 0$, $\eta_t \prec \zeta_t$. Now, since $f_G$ is decreasing and nonnegative,

$$
(5.11) \qquad f_G(\eta_t) \mathbb{1}\{\tau(\eta.) > t\} \geq f_G(\zeta_t) \mathbb{1}\{\tau(\zeta.) > t\}.
$$

Now, we choose $i \neq 0$, and for any $\eta$, we denote by $\zeta := \mathsf{A}_i^+ \eta$. We denote by $E_{\eta,i}^*$ the law of a coupling between $(\eta_\cdot, \zeta_\cdot)$ such that the $\eta$-particles move together. Then,

$$\begin{aligned}
\bar{S}_t^* f_G(\eta) - \bar{S}_t^* f_G(\mathsf{A}_i^+ \eta) &= E_{\eta,i}^*[f_G(\eta_t)\mathbb{1}\{\tau(\eta_\cdot) > t, \tau(\zeta_\cdot) \leq t\} \\
&\quad + (f_G(\eta_t) - f_G(\zeta_t))\mathbb{1}\{\tau(\zeta_\cdot) > t\}] \\
&\leq E_\eta^*[f_G(\eta_t)\mathbb{1}\{\tau(\eta_\cdot) > t\}]P_{0,i}^*(H_0 < \infty) \\
&\quad + E_{\eta,i}^*\left[\left(1 - \frac{f_G(\zeta_t)}{f_G(\eta_t)}\right) f_G(\eta_t)\mathbb{1}\{\tau(\eta_\cdot) > t\}\right].
\end{aligned}$$

We denote by $\gamma_t$ the position at time $t$ of the particle starting in $i$, and thus, $\zeta_t = \mathsf{A}_{\gamma_t}^+ \eta_t$. Now, since $f_G/\psi_\alpha$ is increasing, we have

(5.12) $$\frac{f_G}{\psi_\alpha}(\zeta_t) \geq \frac{f_G}{\psi_\alpha}(\eta_t) \implies \frac{f_G(\zeta_t)}{f_G(\eta_t)} \geq \frac{\psi_\alpha(\mathsf{A}_{\gamma_t}^+ \eta_t)}{\psi_\alpha(\eta_t)} = \frac{\alpha_{\gamma_t}}{\rho}.$$

Thus,

(5.13) $$\begin{aligned}
E_{\eta,i}^*\bigg[\bigg(1 &- \frac{f_G(\zeta_t)}{f_G(\eta_t)}\bigg) f_G(\eta_t)\mathbb{1}\{\tau(\eta_\cdot) > t\}\bigg] \\
&\leq \mathbb{E}_i^*\left[1 - \frac{\alpha_{\gamma_t}}{\rho}\right] E_\eta^*[f_G(\eta_t)\mathbb{1}\{\tau(\eta_\cdot) > t\}].
\end{aligned}$$

Now, we note that, by (1.9),

$$\mathbb{E}_i^*\left[1 - \frac{\alpha_{\gamma_t}}{\rho}\right] = \sum_j P_{0,i}^*(\gamma_t = j)\left(1 - \frac{\alpha_j}{\rho}\right) \leq C_\alpha P_{0,i}^*(H_0 < \infty).$$

Now, for $i \neq 0$, we set

$$\varepsilon_i := 1 \wedge [(1 + C_\alpha) P_{0,i}^*(H_0 < \infty)].$$

For $i = 0$, we set $\varepsilon_i = 1$, and (5.10) holds trivially for any $\eta$. This implies by induction, as in (5.5), that, for any $\underline{\xi} \prec \zeta_t \wedge \zeta_\infty$, and any $\bar{\xi}_t \succ |\zeta_t - \zeta_\infty|$, we have

(5.14) $$|f_t^*(\zeta_t) - f_t^*(\zeta_\infty)| \leq f_t^*(\underline{\xi}) \sum_i \varepsilon_i \bar{\xi}_t(i).$$

Note that, using the arguments in the proof of Lemma 5.1, we deduce from (5.14) that

(5.15) $$\lim_{t \to \infty} \sup_s \left| \int f_s^* \, dT_t(\nu_\rho) - \int f_s^* \, dT_\infty(\nu_\rho) \right| = 0.$$

Indeed, we only need that $\sup_s \|f_s^*\|_{\nu_\rho} < \infty$, which is a simple consequence of (5.10) (see Lemma 7.1 of [2]).

STEP 2. Now, let $(\zeta_t, \zeta_t^\mu)$ be the coupling of Section 4.2 between $T_t(\nu_\rho)$ and $T_t(\mu_G)$. Let $\underline{\xi}_t$ and $\bar{\xi}_t$ be the two independent configurations obtained in Lemma 4.3, and recall that the laws of $\{\underline{\xi}_t, t \geq 0\}$, denoted by $\{\nu_{\beta^t}, t \geq 0\}$, have densities which are uniformly bounded in $L^2(\nu_\rho)$. Then, by taking expectation in (5.14), we obtain, for any $s$ and $t$,

$$(5.16) \qquad \left| \int f_s^* \, dT_t(\nu_\rho) - \int f_s^* \, dT_t(\mu_G) \right| \leq E[f_s^*(\underline{\xi}_t)] \sum_{i \in \mathbb{Z}^d} E[\bar{\xi}_t(i)] \varepsilon_i.$$

Now, by Cauchy–Schwarz,

$$(5.17) \qquad E[f_s^*(\underline{\xi})] \leq \sup_s \|f_s^*\|_{\nu_\rho} \sup_t \left\| \frac{d\nu_{\beta^t}}{d\nu_\rho} \right\|_{\nu_\rho} < \infty.$$

Now, by Corollary 4.4,

$$(5.18) \qquad \lim_{t \to \infty} \sup_s \left| \int f_s^* \, dT_t(\nu_\rho) - \int f_s^* \, dT_t(\mu_G) \right| = 0.$$

By combining (5.15) and (5.18), we conclude that $\{f_t^*, t \geq 0\}$ is an $L^2(\nu_\rho)$-Cauchy sequence with $\lim f_t^* = u^*$.

STEP 3. Let $g$ be as in Proposition 1.4. Arguments similar to those used in Step 2 imply that there is an explicit number $M'(\beta_1)$ such that

$$(5.19) \qquad \left| \int g \, dT_t(\mu_G) - \int g \, dT_t(\nu_\rho) \right| \leq C_g M'(\beta_1) e^{-\beta_1 t} \|g\|_{\nu_\rho},$$

while the proof of Lemma 5.1, with $g$ replacing $u_s$, implies that, for an explicit $M(\beta_1)$,

$$(5.20) \qquad \left| \int g \, dT_t(\nu_\rho) - \int g \, dT_\infty(\nu_\rho) \right| \leq C_g M(\beta_1) e^{-\beta_1 t} \|g\|_{\nu_\rho}.$$

By combining (5.19) and (5.20), we have (1.13) for $\bar{M}(\beta_1) = 2(M(\beta_1) \vee M'(\beta_1))$.

## 6. Bounding the clan.

6.1. *The clan of I-parents.* We bound the number of elements of an $I$-clan, and the clan's width of (3.15), when the particle density $\rho$ is small, and $I = \mathbb{R}$.

By Lemma A.3 of the Appendix, there is $\beta_d$ explicit such that, for any $\beta_1 < \beta_d$, we have a positive density threshold $\rho_c(\beta_1)$ given in (1.6). Henceforth, we consider $\rho \leq \rho_c(\beta_1)$.

PROOF OF LEMMA 3.1. The rectangles $I$-interact only through their time spent on site $\{0\}$, that we have called $\boldsymbol{\sigma}(\mathsf{R})$ in Section 2.1. We can parametrize $\boldsymbol{\sigma}(\mathsf{R})$ by the time of the last visit to $\{0\}$, say, $\bar{s}(\mathsf{R})$, and by $\sigma(\mathsf{R})$, the total length of $\boldsymbol{\sigma}(\mathsf{R})$, whose density $g_\sigma$ is estimated in Lemma A.3. Thus, we will think of the basis of each rectangle to be $(\bar{s}(\mathsf{R}), \sigma(\mathsf{R}))$, rather than the full trajectory $\gamma(\mathsf{R})$. We call this new process the *projected* process $N_p$ defined on $E := \mathbb{R} \times \mathbb{R}^+ \times \mathbb{R} \times \mathbb{R}^+$ equipped with its Borel $\sigma$-field $\mathcal{B}_E$. The elements of $E$ are denoted $x = (s, \sigma, b, l)$. Actually, it is more convenient to reverse time, and think of $s$ as the hitting time of $\{0\}$. We also use, when convenient, $x^s, x^\sigma, x^b$ and $x^l$ for $s, \sigma, b$ and $l$, respectively. We denote by $m$ the support of $N_p$. The density of the intensity measure of $N_p$ is easily seen to be

$$\rho P_{0,0}(H_0 = \infty)\, ds\, g_\sigma(\sigma)\, d\sigma\, db\, \exp(-l)\, dl,$$

in which we will replace henceforth $P_{0,0}(H_0 = \infty)$ by 1, since we only need to bound the clan's size and width.

We consider now $I$-parents, for $I = \mathbb{R}$, of a point $x_0 = (s_0, \sigma_0, b_0, l_0) \in m$. We denote by $\mathbb{P}_1(x_0)$ the domain that $I$-parents of $x_0$ could occupy:

(6.1) $\quad \mathbb{P}_1(x_0) := \{(s, \sigma, b, l) : b < b_0, l > b - b_0, [s, s+\sigma[\,\cap\,[s_0, s_0+\sigma_0[\,\neq \varnothing\}.$

Then, for $A \in \mathcal{B}_E$, we call $M_1(x_0, A)$ the variable counting the parents of $x_0$ falling in $A$:

(6.2) $\quad M_1(x_0, A) := N_p(A \cap \mathbb{P}_1(x_0)) \qquad$ and denote its density by $M_1(x_0, dx)$.

We denote by $m_1(x_0)$ the (random) support of $M_1(x_0, E)$. Now, let $\mathbb{P}_2(x_0)$ be the domain corresponding to the grand-parents of $x_0$:

(6.3) $\quad \mathbb{P}_2(x_0) := \bigcup \{\mathbb{P}_1(x) : x \in m_1(x_0)\} \setminus \mathbb{P}_1(x_0).$

For $A \in \mathcal{B}_E$, we form the counting variable $M_2(x_0, A) := N_p(A \cap \mathbb{P}_2(x_0))$, and denote by $m_2(x_0)$ the support of $M_2(x_0, E)$. We proceed by induction to define, for the $k$-parents of $x_0$,

(6.4) $\quad \mathbb{P}_k(x_0) := \bigcup \{\mathbb{P}_1(x) : x \in m_{k-1}(x_0)\} \setminus \left( \bigcup_{i=1}^{k-1} \mathbb{P}_i(x_0) \right),$

and $M_k(x_0, A) := N_p(A \cap \mathbb{P}_k(x_0))$ with corresponding support $m_k(x_0)$. The $I$-clan of $x_0$ are the points of $\bigcup_k m_k(x_0)$. For $A \in \mathcal{B}_E$, a first obvious bound on $M_k(x_0, A)$ is obtained as we count the parents of the $(k-1)$st generation with their multiplicity:

(6.5) $\quad M_k(x_0, A) \leq \int N_p\left(A \cap \mathbb{P}_1(x) \setminus \left(\bigcup_{i=1}^{k-1} \mathbb{P}_i(x_0)\right)\right) M_{k-1}(x_0, dx),$

and by integrating counting variables over disjoint sets, we have, by independence,

$$E[M_k(x_0, A)] \leq \int E\left[N_p\left(A \cap \mathbb{P}_1(x) \setminus \left(\bigcup_{i=1}^{k-1} \mathbb{P}_i(x_0)\right)\right)\right] E[M_{k-1}(x_0, dx)]$$

$$\leq \int E[N_p(A \cap \mathbb{P}_1(x))] E[M_{k-1}(x_0, dx)],$$

where the second inequality corresponds to counting the parents of points of $m_{k-1}(x_0)$ even if they are part of an earlier generation. By induction, we obtain the upper bound

(6.6)
$$E[M_k(x_0, A)] \leq \int \cdots \int E[M_1(x_{k-1}, A)]$$
$$\times E[M_1(x_{k-2}, dx_{k-1})] \cdots E[M_1(x_0, dx_1)].$$

The following estimates provide an upper bound for (6.6): consider the cylinder $dA_\sigma := \mathbb{R} \times [\sigma, \sigma + d\sigma[ \times \mathbb{R} \times \mathbb{R}^+$, infinitesimal in the $\sigma$-direction, and

(6.7)
$$E[M_1(x_0, dA_\sigma)]$$
$$= \int \cdots \int \mathbb{1}\{b < b_0, l > b_0 - b, [s, s + \sigma[ \cap [s_0, s_0 + \sigma_0[ \neq \varnothing\}$$
$$\times \rho\, ds\, e^{-l}\, db\, dl\, g_\sigma(\sigma)\, d\sigma$$
$$= \rho \int_{-\infty}^{\infty} \mathbb{1}\{[s, s + \sigma[ \cap [s_0, s_0 + \sigma_0[ \neq \varnothing\}\, ds\, g_\sigma(\sigma)\, d\sigma$$
$$= \rho(\sigma_0 + \sigma)g_\sigma(\sigma)\, d\sigma.$$

The expression (6.7) depends only on $\sigma_0$ and $\sigma$. Thus, when performing the integral of (6.6), we first integrate $x_i^s, x_i^b$ and $x_i^l$, for $i = 1, \ldots, k$, over $\mathbb{R} \times \mathbb{R} \times \mathbb{R}^+$. If we call $f(\sigma', \sigma) = \rho(\sigma' + \sigma)g_\sigma(\sigma)$, we have from (6.6) a bound on the number of $\mathbb{R}$-parents of the $k$th generation,

(6.8) $$E[M_k(x_0, E)] \leq \int \cdots \int f(\sigma_0, \sigma_1)f(\sigma_1, \sigma_2) \cdots f(\sigma_{k-1}, \sigma_k)\, d\sigma_1 \cdots d\sigma_k.$$

Note that the following simple inequality,

(6.9)
$$(\sigma_0 + \sigma_1)(\sigma_1 + \sigma_2) \cdots (\sigma_{k-1} + \sigma_k)$$
$$\leq (1 + \sigma_0 + \sigma_0^2) \cdots (1 + \sigma_k + \sigma_k^2),$$

implies that $E[M_k(x_0, E)] \leq (\rho E[1 + \sigma + \sigma^2])^k (1 + \sigma_0 + \sigma_0^2)$. Since $E[\exp(\beta_1 \times \sigma)] < \infty$, by Lemma A.3, we have $E[1 + \sigma + \sigma^2] < \infty$, and for $\rho < \rho_c(\beta_1)$, the $I$-clan is $\mathbb{Q}$-a.s. finite. $\square$

We define the width of the clan of $x_0$ for the *projected* point process $N_p$ [compare with (3.15)] by

(6.10) $$W(x_0) = \bigcup \left\{ [x^s, x^s + x^\sigma] : x \in \bigcup_k m_k(x_0) \right\}.$$

Note that in Section 4.1 we have introduced $N_p$ with $s$ constrained in some time interval $[-t, 0]$ rather than $\mathbb{R}$. The following lemma is similar to [8].

LEMMA 6.1. *For $\beta_1 < \beta_d$ of Lemma A.1, and $\rho < \rho_c(\beta_1)$, we have, for any $t \in \mathbb{R}$,*

(6.11) 
$$\forall s_0, b_0, l_0 \in \mathbb{R} \times \mathbb{R} \times \mathbb{R}^+,$$
$$\int P(t \in W((s_0, \sigma_0, b_0, l_0))) g_\sigma(\sigma_0) \, d\sigma_0 \leq \frac{e^{-\beta_1 |t-s_0|}}{\rho_c(\beta_1) - \rho}.$$

PROOF. Looking at the definition (6.10), it is clear that if the width of the clan contains $t$, then at some generation the total width contains $t$. In other words,

$$\mathbb{1}\{t \in W(x_0)\} \leq \sum_{k \geq 0} \mathbb{1}\{t \in [x^s, x^s + x^\sigma[ : x \in m_k(x_0)\}$$

[where we set $m_0(x_0) := \{x_0\}$]

(6.12)
$$\leq \sum_{k \geq 0} \int \cdots \int \mathbb{1}\{x_0^\sigma + \cdots + x_k^\sigma > |t - x_0^s|\}$$
$$\times \prod_{i=1}^k M_1(x_{i-1}, dx_i).$$

To compute the expectation of the right-hand side of (6.12), we first integrate $x_i^s, x_i^b, x_i^l$ for $i = 1, \ldots, k$. For $x_0 = (s_0, \sigma_0, b_0, l_0)$, we obtain, using (6.8), that

(6.13)
$$P(t \in W(x_0)) \leq \sum_{k \geq 1} \int \cdots \int \mathbb{1}\{\sigma_0 + \cdots + \sigma_k > |t - s_0|\}$$
$$\times \prod_{i=1}^k (f(\sigma_{i-1}, \sigma_i) \, d\sigma_i) + \mathbb{1}\{\sigma_0 > |t - s_0|\}.$$

If we define

(6.14) $$\psi(t) := \int P(t \in W(x_0)) g_\sigma(\sigma_0) \, d\sigma_0,$$

then, from (6.13), we obtain

$$
(6.15) \quad \begin{aligned} \psi(t) &\leq e^{-\beta_1|t-s_0|} \\ &\times \left( \sum_{k\geq 0} \int \cdots \int e^{\beta_1(\sigma_0+\cdots+\sigma_k)} g_\sigma(\sigma_0) \prod_{i=1}^{k} (f(\sigma_{i-1},\sigma_i)\,d\sigma_i)\,d\sigma_0 \right). \end{aligned}
$$

Using (6.9), we obtain, for $\rho < \rho_c(\beta_1)$,

$$
(6.16) \quad \psi(t) \leq e^{-\beta_1|t-s_0|} \left( \frac{1}{\rho_c(\beta_1)} + \sum_{k\geq 1} \left(\frac{\rho}{\rho_c(\beta_1)}\right)^k \right) \leq \frac{e^{-\beta_1|t-s_0|}}{\rho_c(\beta_1) - \rho}. \quad \square
$$

6.2. *Bad parents.* We consider the point process $\tilde{N}$ introduced in Section 4.2, on $E^+ := \mathbb{R}^+ \times \mathbb{R}^+ \times \mathbb{R} \times \mathbb{R}^+$, whose density of the intensity measure is given in (4.27). In this section we evaluate the probability that a point $x_0 = (s_0, \sigma_0, b_0, l_0)$ has a bad parent. In other words, we estimate the event

$$
(6.17) \quad K(x_0) = \left\{ \tilde{N}\left( \left\{ \bigcup \mathbb{P}_1(x) : x \in m_k(x_0), k \in \mathbb{N} \right\} \right) \geq 1 \right\}.
$$

LEMMA 6.2. *For $\beta_1 < \beta_d$ of Lemma A.1, and $\rho < \rho_c(\beta_1)$, we have*

$$
(6.18) \quad \begin{aligned} \forall\, x_0 &= (s_0, \sigma_0, b_0, l_0) \in \mathbb{R}^+ \times \mathbb{R}^+ \times \mathbb{R} \times \mathbb{R}^+, \\ P(K(x_0)) &\leq \frac{\rho c_0(\beta_1)}{\beta_1(1 - \rho/\rho_c(\beta_1))} e^{-\beta_1 s_0}. \end{aligned}
$$

PROOF. First note that [with $m_0(x_0) = \{x_0\}$]

$$
(6.19) \quad \begin{aligned} \mathbb{1}\{K(x_0)\} &\leq \sum_{k=0}^{\infty} \tilde{N}\left( \left\{ \bigcup \mathbb{P}_1(x) : x \in m_k(x_0) \right\} \right) \\ &\leq \sum_{k=0}^{\infty} \int \tilde{N}(\mathbb{P}_1(x)) M_k(x_0, dx). \end{aligned}
$$

Thus, using independence of $N_p$ and $\tilde{N}$, and the bound (6.6),

$$
(6.20) \quad P(K(x_0)) \leq \sum_{k=0}^{\infty} \int E[\tilde{N}(\mathbb{P}_1(x_k))] E[M_1(x_{k-1}, dx_k)] \cdots E[M_1(x_0, dx_1)].
$$

We first integrate $\tilde{N}$ over $dA_\sigma \cap \mathbb{P}_1(x_k)$ to obtain

$$
\begin{aligned} E[\tilde{N}(dA_\sigma \cap \mathbb{P}_1(x_k))] &\leq \int_0^\infty ds\, \rho c_0(\beta_1) \exp(-\beta_1 s) \\ &\quad \times \mathbb{1}\{[s, s+\sigma[\,\cap\,[s_k, s_k+\sigma_k[\,\neq \varnothing\} g_\sigma(\sigma)\,d\sigma \end{aligned}
$$

$$\text{(6.21)} \qquad \leq \rho c_0(\beta_1) \int_{(s_k-\sigma)^+}^{(s_k+\sigma_k)^+} ds \exp(-\beta_1 s) g_\sigma(\sigma) \, d\sigma$$

$$\leq \frac{c_0(\beta_1)\rho}{\beta_1} \exp(-\beta_1(s_k-\sigma)^+) g_\sigma(\sigma) \, d\sigma.$$

Note that this is independent of $b_k, l_k$ and $\sigma_k$. Thus, after integrating the bad points intensity,

$$\text{(6.22)} \qquad P(K(x_0)) \leq \sum_{k=0}^{\infty} \frac{c_0(\beta_1)\rho}{\beta_1} \int \cdots \int_0^\infty \exp(-\beta_1(s_k-\sigma)^+) g_\sigma(\sigma) \, d\sigma$$

$$\times \prod_{i=1}^{k} E[M_1(x_{i-1}, dx_i)].$$

Note that $|s_0 - s_k| \leq \sigma_1 + \cdots + \sigma_k$, so that with the notation of the proof of Lemma 6.1,

$$P(K(x_0)) \leq \sum_{k=1}^{\infty} \frac{c_0(\beta_1)\rho}{\beta_1} \int \cdots \int_0^\infty \exp\left(-\beta_1\left(s_0 - \sum_{i=1}^{k} \sigma_i\right)^+\right)$$

$$\times \prod_{i=1}^{k} g_\sigma(\sigma_i) f(\sigma_{i-1}, \sigma_i) \, d\sigma_i$$

$$\text{(6.23)} \qquad \leq \sum_{k=1}^{\infty} \frac{c_0(\beta_1)\rho}{\beta_1} \int \cdots \int_0^\infty \exp\left(-\beta_1\left(s_0 - \sum_{i=1}^{k} \sigma_i\right)\right)$$

$$\times \prod_{i=1}^{k} (\rho g_\sigma(\sigma_i)(1 + \sigma_i + \sigma_i^2) \, d\sigma_i)$$

$$\leq e^{-\beta_1 s_0} \sum_{k=1}^{\infty} \frac{c_0(\beta_1)\rho}{\beta_1} \left(\frac{\rho}{\rho_c(\beta_1)}\right)^k$$

$$= \frac{c_0(\beta_1)\rho}{\beta_1(1 - \rho/\rho_c(\beta_1))} e^{-\beta_1 s_0}. \qquad \square$$

6.3. *The clan of $(I, \mu)$-parents.* In this section we prove that $\mathbf{B}^{R,U}(I, \mu)$ is a.s. finite.

PROOF OF LEMMA 3.8. If $R$ is the range of the Gibbs measure $\mu_G$, we define $K := (2R+1)^d$. We choose $I = [-t, 0]$, and for ease of notation, we set $s = -t$.

As in [7], we consider at time $\mathfrak{b} = 0$ one rectangle, $\mathsf{R}_0$, and build its clan backward in time. For simplicity, we work with positive backward time. For each (backward) time $\mathfrak{b}$, we build a set of rectangles, denoted by $\mathbb{B}_\mathfrak{b}$ with the property that if at a certain time $\mathfrak{b}$ we have $\mathbb{B}_\mathfrak{b} = \varnothing$, then none of the parents of $\mathsf{R}_0$ are alive at time $\mathfrak{b}$. If we denote by $\tau_\varnothing$ the first time where $\mathbb{B}_\mathfrak{b} = \varnothing$, then we show that $E[\tau_\varnothing] < \infty$.

First, $\mathbb{B}_0 := \{\mathsf{R}_0\}$. Then, for each small $\delta > 0$, we partition $\mathbb{Z}^d$ into $\mathbb{D}_\delta := \{i \in \mathbb{Z}^d : \alpha_i/\rho \geq 1 - \delta\}$, and its complement $\mathbb{D}_\delta^c$. Note that $\mathbb{D}_\delta^c$ is bounded since $\sum(1 - \alpha_i/\rho)^2 < \infty$, by Remark 1.5. The point of this partition of $\mathbb{Z}^d$ is that a rectangle $\mathsf{R}$ with $\gamma_s(\mathsf{R}) \in \mathbb{D}_\delta$ have a small probability [less than $\delta$ by (3.23)] to have a $\mu$-parent.

Now, $\mathbb{B}_\mathfrak{b}$ contains all rectangles of $\mathbf{C}$ whose trajectory is in $\mathbb{D}_\delta^c$ at time $s$ and whose life-epoch contains $\mathfrak{b}$. Thus, it is convenient to associate with each site $i$ of $\mathbb{D}_\delta^c$ a birth and death process of intensity $\rho$, and to attach, to each birth-time, a trajectory drawn from $dP_{s,i}(\gamma)$, which we color in blue. The trajectories of rectangles in $\mathbb{B}_\mathfrak{b}$ with a position in $\mathbb{D}_\delta$ at time $s$ are colored in yellow.

It is also convenient to think that rectangles of $\mathbb{B}$ generate $(I, \mu)$-parents at their death time. This does not lengthen the life-time of parents, since the exponential life-time $\tau$ satisfies $P(\tau > t + s | \tau > s) = P(\tau > t)$, but we underestimate the number of parents alive at a given time. However, since we are only interested in showing that there is a finite time at which the clan dies out, the life-times we are ignoring are insignificant since their children are alive at the overlapping times.

We recall that to build the $I$-parents, one only considers $\{\sigma(\mathsf{R}), \mathsf{R} \in \mathbf{C}\}$. Now, in the stationary rectangle process, a rectangle $\mathsf{R}_0$ has a Poisson number of parents with $\sigma(\mathsf{R}) \in [\sigma, \sigma + d\sigma[$ with an intensity measure whose density $m(\mathsf{R}_0, \sigma)$ is bounded by $\rho(\sigma(\mathsf{R}_0) + \sigma)g_\sigma(\sigma)$. Note that this bound only depends on $\sigma(\mathsf{R}_0)$, and that the distribution of $\sigma(\mathsf{R}_0)$ (once we assume the trajectory has touched $\{0\}$) is independent of $\gamma(\mathsf{R}_0)_s$. Thus, the only relevant properties of a trajectory are its time-width $\sigma(\mathsf{R})$, and its location at time $s$ (actually only whether it is blue or yellow). Now, we overcount the number of $I$-parents when we assume that each point has an independent Poisson number of parents, all of them colored yellow. Indeed, we do not need to worry about the blue ones, since we have included them all in $\{\mathbb{B}_\mathfrak{b}, \mathfrak{b} \geq 0\}$. To make things easier, we actually discretize the possible values of the time-width. Thus, a rectangle $\mathsf{R}_0$ with $\sigma(\mathsf{R}_0) \in [k-1, k[$ gives rise to a Poisson number of $I$-parents with $\sigma(\mathsf{R}) \in [i-1, i[$ with intensity measure bounded by

$$m(k, i) \leq \rho(k + i)q_i \qquad \text{with } q_i := \int_{i-1}^{i} g_\sigma(\sigma)\, d\sigma.$$

We can simplify the description of the above-mentioned birth and death process giving rise to the blue trajectory. We actually consider, at each site of $\mathbb{D}_\rho^c$, a Poisson process of intensity $\rho$ and we associate with every mark a time-width variable with distribution $\{q_i\}$. This procedure has the effect of overestimating the parents number, since a trajectory can very well not touch $\{0\}$ during the time-period $I$. The configuration of blue marks is denoted by $\beta : 1, 2, \ldots \to \mathbb{N}$, where $\beta(i)$ is the number of blue marks with a time-width in $[i-1, i[$. Similarly, the configuration of yellow marks is denoted by $y$.

We describe now the full evolution of a rectangle at its death-time:

- If it is blue, we assume it gives rise to $K$ yellow points with an independent distribution of the time-width drawn from $\{q_i, i \geq 1\}$. It also gives rise to $I$-parents as described above.
- If it is yellow, with probability $\delta$, it behaves as a blue point, and with probability, $1 - \delta$ it has only $I$-parents.

Thus, we are giving to trajectories touching $\mathbb{D}_\delta$ at time $s$ more parents than what comes from the prescription of detailed balance. The advantage is that we do not keep track of the whole trajectory, but only of its color.

We write now the generator of the evolution of colored rectangles in set $\mathbb{B}$ backward in time. The configuration variable is $x = (\beta, y)$ with $\beta, y \in \mathbb{N}^{\{1,2,\ldots\}}$. We denote by $\underline{i} = (i_1, \ldots, i_K)$, where $i_j \in \{1, 2, \ldots\}$, and we use $\mathsf{A}_i^+$ (resp. $\mathsf{A}_i^-$) for the action of adding (resp. canceling) a mark with time-width in $[i-1, i[$. For a function $f$ of $(\beta, y)$,

(6.24) $$\mathcal{L}f(\beta, y) = \bar{\mathcal{L}}f(\beta, y) + \underline{\mathcal{L}}f(\beta, y),$$

where $\bar{\mathcal{L}}$ accounts for the evolution of blue parents, and $\underline{\mathcal{L}}$ accounts for the yellow parents,

(6.25) $$\begin{aligned}\bar{\mathcal{L}}f(\beta, y) &= \rho|\mathbb{D}_\delta^c| \sum_{i \geq 1} q_i(f(\mathsf{A}_i^+ \beta, y) - f(\beta, y)) \\ &+ \sum_{k \geq 1} \beta(k) \sum_{\underline{i}} \sum_{\zeta : \sum_j \zeta(j) < \infty} \left(\prod_{j=1}^K q_{i_j}\right) Q(k, \zeta) \\ &\quad \times \left(f\left(\mathsf{A}_k^- \beta, \prod_{j=1}^K \mathsf{A}_{i_j}^+ y + \zeta\right) - f(\beta, y)\right),\end{aligned}$$

where, for a configuration of parents $\zeta$, we set $Q(k, \zeta) = \prod_{j \in \mathbb{N}} e^{-m(k,j)} m(k,j)^{\zeta(j)}/\zeta(j)!$, and

$$\underline{\mathcal{L}}f(\beta, y) = \sum_{k \geq 1} y(k)\delta \sum_{\underline{i}} \sum_{\zeta : \sum_j \zeta(j) < \infty} \left(\prod_{j=1}^K q_{i_j}\right) Q(k, \zeta)$$

$$\text{(6.26)} \qquad \times \left( f\left(\beta, \left(\prod_{j=1}^{K} \mathsf{A}_{i_j}^+\right) \mathsf{A}_k^- y + \zeta\right) - f(\beta, y) \right)$$

$$+ \sum_{k \geq 1} y(k)(1-\delta)$$

$$\times \sum_{i} \sum_{\zeta : \sum_j \zeta(j) < \infty} Q(k, \zeta)(f(\beta, \mathsf{A}_k^- y + \zeta) - f(\beta, y)).$$

Now, we look for a Lyapounov function, following the classical Foster's arguments. We consider the function $f(\beta, y) = \sum_j \varphi_j (C\beta(j) + y(j))$, with $\varphi_j = \sqrt{1 + j + j^2}$ and $C$ a positive (large) constant to be chosen later. With this choice of $f$, simple algebra yields

$$\bar{\mathcal{L}} f(\beta, y) = C\rho |\mathbb{D}_\delta^c| \sum_{i \geq 1} q_i \varphi_i - C \sum_k \beta(k) \varphi_k$$

$$\text{(6.27)} \qquad + K \left( \sum_k \beta(k) \right) \left( \sum_i q_i \varphi_i \right) + \sum_{k,i} \beta(k) m(k,i) \varphi_i,$$

and

$$\text{(6.28)} \quad \underline{\mathcal{L}} f(\beta, y) = -\sum_k y(k) \varphi_k + K\delta \left( \sum_k y(k) \right) \sum_{i \geq 1} q_i \varphi_i + \sum_{k,i} y(k) m(k,i) \varphi_i.$$

Now, $m(k,i) \leq \rho(k+i) q_i \leq \rho \varphi_k \varphi_i q_i$. By Lemma A.3, $c_0 := \sum_i q_i \varphi_i^2 < \infty$, and we obtain

$$\sum_{k,i} y(k) m(k,i) \varphi_i \leq \rho c_0 \sum_k y(k) \varphi_k \quad \text{and}$$

$$\text{(6.29)}$$

$$\sum_{k,i} \beta(k) m(k,i) \varphi_i \leq \rho c_0 \sum_k \beta(k) \varphi_k.$$

Moreover, for $\rho < \rho_c(\beta_1)$, using that $\varphi_i \leq 2\varphi_{i-1}$,

$$\sum_{i \geq 1} q_i \varphi_i \leq c_1 := 2(E[1 + \sigma + \sigma^2])^{1/2} < \infty \quad \text{and}$$

$$\text{(6.30)}$$

$$\sum_k y(k) \leq \sum_k y(k) \varphi_k.$$

Thus,

$$\text{(6.31)} \qquad \underline{\mathcal{L}} f(\beta, y) \leq -(1 - \rho c_0 - \delta K c_1) \sum_k y(k) \varphi_k.$$

Also, with similar computations,

$$\text{(6.32)} \qquad \bar{\mathcal{L}} f(\beta, y) \leq C\rho |\mathbb{D}_\delta^c| c_1 - (C - K c_1 - \rho c_0) \sum_k \beta(k) \varphi_k.$$

First, we choose $\rho < \rho_c(\beta_1)$ so that $\rho c_0 < 1$ (this can always be achieved by making the discretization fine enough). Second, we choose $\delta$ so that $c_3 := 1 - \rho c_0 - \delta K c_1 > 0$, and $C > 1$ such that $c_4 := C - Kc_1 - \rho c_0 > 0$. Then, if we set $c_5 := C\rho|\mathbb{D}_\delta^c|c_1$, we obtain

$$(6.33) \quad \mathcal{L}f(\beta,y) \leq c_5 - c_4 \sum_k \beta(k)\varphi_k - c_3 \sum_k y(k)\varphi_k \leq c_5 - c_6 f(\beta,y),$$

where $c_6 = \min(c_3, c_4/C)$. Let $k_0$ be such that $c_6 \varphi_{k_0} - c_5 > 0$, and define

$$Z = \{(\beta,y): y(k) = \beta(k) = 0 \text{ for } k > k_0 \text{ and } y(i) \vee C\beta(i) < M_i, i = 1,\ldots,k_0\},$$

where $M_i$ are such that $c_6 \min\{M_i\varphi_i : i \leq k_0\} - c_5 > 0$. Note that $Z$ has finitely many configurations. Define

$$c_7 := c_6 \min(\varphi_{k_0}, \min(\{M_i\varphi_i : i \leq k_0\})) - c_5 > 0,$$

and $\tau_Z = \inf\{s : x_s \in Z\}$, and note that on $\{\tau_Z > t\}$, we have from (6.33) that $\mathcal{L}f(x_t) \leq -c_7$. Now, for any state $x = (\beta,y)$, we consider the mean-zero martingale

$$(6.34) \quad \begin{aligned} M_t &= f(x_{t\wedge\tau_Z}) - f(x) - \int_0^t \mathcal{L}f(x_s)\mathbb{1}\{\tau_Z > s\}\,ds \\ &\geq f(x_{t\wedge\tau_Z}) - f(x) + c_7(t\wedge\tau_Z). \end{aligned}$$

We take the expectation of each side of (6.34) and take $t$ to infinity to obtain $f(x) \geq c_7 E_x[\tau_Z]$. Now, for each $z \in Z$, the probability of reaching the empty configuration $\{y \equiv 0, \beta \equiv 0\}$ in a unit-time interval is positive. Finally, a standard renewal argument yields that $E_x[\tau_\varnothing] < \infty$.

## APPENDIX

To ease the reading, we first derive some classical bound for $P_{0,0}(t < H_0 < \infty)$. We denote by $\{S_n, n \in \mathbb{N}\}$ the discrete sums $S_n = \gamma_1 + \cdots + \gamma_n$, where the $\gamma_i$ are i.i.d. with law $\{p(0,\cdot)\}$. We introduce some definitions with the notation of [12]. For $z \in \mathbb{R}^d$, the finite range assumption on $\{p(0,i), i \in \mathbb{Z}^d\}$ implies that the exponential moments of the increments of $\gamma_1$ exist and

$$\forall z \in \mathbb{R}^d \quad \Phi(z) := E[e^{z\cdot\gamma_1}] = \sum_{i\in\mathbb{Z}^d} p(0,i)e^{z\cdot i} \quad \text{and}$$

(A.1)
$$D := \{z \in \mathbb{R}^d : \Phi(z) \leq 1\}.$$

The finite range and irreducibility assumptions imply that $\Phi$ is well defined and strictly positive on $\mathbb{R}^d$. It is shown in [9] (see also [12], Lemma 1.1) that $D$ is compact and convex, that $\nabla\Phi$ does not vanish on $\partial D := \{z : \Phi(z) = 1\}$, and that $z \mapsto \nabla\Phi(z)/\|\nabla\Phi(z)\|$ is a continuous bijection from $\partial D$ to the unit sphere of $\mathbb{R}^d$. A simple consequence is that $D \setminus \partial D$ is not empty. Indeed, by contradiction, assume that $D = \partial D$ and let $z^* \neq 0 \in \partial D$. Then, for any $z \in D$ and $t \in [0,1]$,
$$\Phi(tz) = 1 \quad \text{and} \quad \Phi(z + t(z^* - z)) = 1.$$

Thus, by differentiating, for any $t, s \in [0,1]$,

(A.2) $\qquad \nabla \Phi(tz) \cdot z = 0 \quad \text{and} \quad \nabla \Phi(z + s(z^* - z)) \cdot (z^* - z) = 0.$

We choose $t = 1$ and $s = 0$ in (A.2), and add the two gradients to obtain
$$\nabla \Phi(z) \cdot z^* = 0 \qquad \forall z \in \partial D.$$

This contradicts that $\nabla \Phi(\cdot)/\|\nabla \Phi\| : \partial D \to S^{d-1}$ is bijective. Thus, there is $z_0 \in D$ such that

(A.3) $\qquad 0 < \Phi(z_0) = \inf\{\Phi(z)\} < 1.$

We denote by $\tilde H_0 = \inf\{n \geq 1 : S_n = 0\}$ and by $H_0$ the analogue for continuous-time walks. In other words, $H_0 = \infty$ if $\tilde H_0 = \infty$, and otherwise,

(A.4) $\qquad H_0 = \sum_{i=1}^{\tilde H_0} \tau_i,$

where $\{\tau_i, i \in \mathbb{N}\}$ are i.i.d. exponential times of intensity 1. From (A.3), we obtain the following estimates.

LEMMA A.1. *Let $\beta = 1 - \Phi(z_0)$ (with $0 < \beta < 1$). Then,*

(A.5)
$$\int [e^{\beta H_0(\gamma)} \mathbb{1}\{H_0(\gamma) < \infty\}] \, dP_{0,0}(\gamma) \leq 1 \quad \text{and}$$
$$P_{0,0}(t < H_0 < \infty) \leq e^{-\beta t}.$$

PROOF. We first work in discrete time. Form the martingale $M_n = \exp(z_0 \cdot S_n)/\Phi(z_0)^n$. Note that $\{M_{\tilde H_0 \wedge n}, n \in \mathbb{N}\}$ is a positive martingale [though $P_{0,0}(\tilde H_0 = \infty) > 0$]. Thus,

(A.6) $\quad 1 = E_{0,0}[M_{\tilde H_0 \wedge n}] = \dfrac{E_{0,0}[e^{z_0 \cdot S_n} \mathbb{1}\{\tilde H_0 > n\}]}{\Phi(z_0)^n} + E_{0,0}\left[\dfrac{\mathbb{1}\{\tilde H_0 \leq n\}}{\Phi(z_0)^{\tilde H_0}}\right].$

As we take the limit $n$ to infinity in (A.6), we obtain

(A.7) $\qquad 1 \geq E_{0,0}\left[\dfrac{\mathbb{1}\{\tilde H_0 < \infty\}}{\Phi(z_0)^{\tilde H_0}}\right].$

Note that this implies, by Chebyshev's inequality, $P_{0,0}(n < \tilde H_0 < \infty) \leq \Phi(z_0)^n$.

Using (A.4) and (A.7), we obtain

(A.8)
$$E_{0,0}[\exp(\beta H_0)\mathbb{1}\{H_0 < \infty\}] = E_{0,0}\left[\left(\frac{1}{1-\beta}\right)^{\tilde{H}_0}\mathbb{1}\{\tilde{H}_0 < \infty\}\right]$$
$$= E_{0,0}\left[\frac{\mathbb{1}\{\tilde{H}_0 < \infty\}}{\Phi(z_0)^{\tilde{H}_0}}\right] \leq 1.$$

The second inequality in (A.5) is a direct consequence of Chebyshev's inequality. □

We need now a bound on the width of a point, which we had called $\sigma$. We decompose a walk starting at 0 into its renewal parts:

- Let $\{Y^{(i)}, i \in \mathbb{N}\}$ be i.i.d. with law $\{p(0, \cdot)\}$ representing the first random move away from 0.
- Let $\{\tau^{(i)}, i \in \mathbb{N}\}$ be i.i.d. exponential times of mean 1, representing the waiting times at 0 (before doing the move $Y^{(i)}$).
- Let $\{\{\gamma_s^{(i)}, s \geq 0\}, i \in \mathbb{N}\}$ be independent walks with transition $\{p(i,j)\}$ and $\gamma_0^{(i)} = 0$.

The convex hull of $\Sigma(\gamma)$ is made up by adding the successive excursions times and waiting times in 0. Let $H^{(i)}$ be the $i$th excursion time, $H^{(i)} = \inf\{s > 0 : \gamma_s^{(i)} + Y^{(i)} = 0\}$, and denote the label of the last excursion by $\kappa = \sup\{i \in \mathbb{N} : H^{(i)} < \infty\}$. Note that $\kappa$ is a geometric variable with $P(\kappa = n) = P_{0,0}(H_0 < \infty)^n P_{0,0}(H_0 = \infty)$. Then,

$$\sigma = \mathbb{1}\{\kappa = 0\}(\tau^{(0)}) + \mathbb{1}\{\kappa = 1\}(\tau^{(0)} + H^{(1)} + \tau^{(1)}) + \cdots$$
$$+ \mathbb{1}\{\kappa = n\}\left(\tau^{(0)} + \sum_{i=1}^{n}(H^{(i)} + \tau^{(i)})\right) + \cdots.$$

Note that since $P_{0,0}(H_0 < \infty) < 1$, we have (with $\tau$ denoting an exponential time)

$$G_\sigma(z) := E_{0,0}[e^{z\sigma}] = E_{0,0}[e^{z\tau}]\sum_{i \geq 0} P(\kappa = i)(E_{0,0}[e^{zH_0}\mathbb{1}\{H_0 < \infty\}]E[e^{z\tau}])^i$$

$$= \frac{P_{0,0}(H_0 = \infty)}{1-z}\frac{1}{1 - P_{0,0}(H_0 < \infty)E_{0,0}[e^{zH_0}\mathbb{1}\{H_0 < \infty\}]E[e^{z\tau}]}.$$

Thus, as a simple consequence of Lemma A.1, we have the following estimate.

LEMMA A.3. *Let $\beta$ be as in Lemma* A.1. *If we define the positive constant*

$$\beta_d := \min(\beta, P_{0,0}(H_0 = \infty)),$$

(A.9)
$$\text{then, for } z < \beta_d, \ G_\sigma(z) := E_{0,0}[\exp(z\sigma)] < \infty.$$

To control discrepencies, we need the following simple estimates.

LEMMA A.4. *Let $d \geq 3$. Then:*

(i) $\sum_{i \neq 0} \sum_j P_{0,i}(H_0 < \infty) p(0,j) P_{0,j}(\gamma_t = i, H_0 > t) \leq e^{-\beta_d t}$.

(ii) $\sum_{i \neq 0} P_{0,i}(H_0 < \infty) P^*_{0,i}(t < H_0 < \infty) \leq e^{-\beta_d t}$.

(iii) $\lim_{t \to \infty} \sum_{i \neq 0} \sum_j P^*_{0,i}(H_0 < \infty) p(0,j) P_{0,j}(\gamma_t = i, H_0 > t) = 0$.

(iv)

(A.10)
$$\lim_{t \to \infty} \sum_{i \neq 0} P_{0,i}(H_0 < \infty) P_{0,i}(t < H_0 < \infty) = 0.$$

(v) *For any $\beta_1 < \beta_d$, there is $c_1(\beta_1)$ such that*
$$\sum_{i \in \mathbb{Z}^d} P_{0,0}(\gamma_t = i) P_{0,i}(H_0 < \infty) \leq c_1(\beta_1) e^{-\beta_1 t}.$$

*Moreover, if $\{\varepsilon_i\}$ is such that $\sum_i \varepsilon_i^2 < \infty$, and if we replace $P_{0,i}(H_0 < \infty)$ or $P^*_{0,i}(H_0 < \infty)$ in (i)–(iv) by $\varepsilon_i$, then the limit as $t$ tends to infinity is zero.*

PROOF. (i) We first show the discrete version of (i). We fix an integer $n$, condition on $S_n$ using the Markov property and Lemma A.1,

(A.11)
$$\Phi(z_0)^n \geq P_{0,0}(n < \tilde{H}_0 < \infty)$$
$$= \sum_{i \neq 0} P_{0,0}(S_n = i, \tilde{H}_0 > n) P_{0,i}(\tilde{H}_0 < \infty).$$

To pass to continuous time, let $N_t$ be the Poisson number of marks before $t$, and decompose over the possible values of $N_t$, for $j$ with $p(0,j) \neq 0$,

$$P_{0,j}(\gamma_t = i, H_0 > t) = \sum_{n \in \mathbb{N}} P_{0,j}(N_t = n, S_n = i, \tilde{H}_0 > n)$$
$$= \sum_{n=0}^{\infty} P(N_t = n) P_{0,j}(S_n = i, \tilde{H}_0 > n).$$

Thus,
$$\sum_{j, i \neq 0} P_{0,i}(H_0 < \infty) p(0,j) P_{0,j}(\gamma_t = i, H_0 > t)$$

$$= \sum_{n\geq 1,i} P(N_t = n-1)P_{0,0}(S_n = i, \tilde{H}_0 > n)P_{0,i}(\tilde{H}_0 < \infty)$$

$$= \sum_{n=1}^{\infty} P(N_t = n-1)P_{0,0}(n < \tilde{H}_0 < \infty) \leq \Phi(z_0)E_{0,0}[\Phi(z_0)^{N_t}]$$

$$\leq \exp(-t(1-\Phi(z_0))).$$

(ii) Similarly, it is enough to prove the discrete version of (ii). By reversing time, note that, for $i \neq 0$, $P^*_{0,i}(\tilde{H}_0 = n) = P_{0,0}(S_n = i, \tilde{H}_0 > n)$. Thus, by (A.11),

(A.12) $$\sum_{i\neq 0} P_{0,i}(\tilde{H}_0 < \infty)P^*_{0,i}(\tilde{H}_0 = n) = P_{0,0}(n < \tilde{H}_0 < \infty).$$

The extention to continuous time is done as in point (i).

(iii) By reversing time,

(A.13) $$\sum_{i\neq 0}\sum_j P^*_{0,i}(H_0 < \infty)p(0,j)P_{0,j}(S_n = i, \tilde{H}_0 > n) \leq \sum_{i\neq 0} f_n(i),$$

with $f_n(i) := P^*_{0,i}(H_0 < \infty)P^*_{0,i}(\tilde{H}_0 = n+1)$. Note that, for any fixed $i$, $f_n(i)$ tends to 0 as $n$ tends to infinity. Now, since

(A.14) $\quad f_n(i) \leq P^*_{0,i}(H_0 < \infty) \quad \text{and in } d \geq 3, \quad \sum (P^*_{0,i}(H_0 < \infty))^2 < \infty,$

Lebesgue dominated convergence yields the discrete version (iii). The passage to continuous time is similar to (i). Point (iv) presents now no difficulty. We omit its proof.

(v) First note that $\gamma_t = 0$ implies that $\sigma \geq t$, thus,

$$\sum_{i\in\mathbb{Z}^d} P_{0,0}(\gamma_t = i)P_{0,i}(H_0 < \infty) \leq P_{0,0}(\sigma \geq t) + \sum_{i\neq 0} P_{0,0}(\gamma_t = i)P_{0,i}(\tilde{H}_0 < \infty).$$

Note that, by Lemma A.3, $P_{0,0}(\sigma \geq t) \leq G_\sigma(\beta_d)\exp(-\beta_d t)$. Now, we deal with the discrete walk, and show that, for any $\delta_1$ with $1 < \delta_1 < 1/\Phi(z_0)$, there is a number $M$ such that, for any integer $n$,

$$\sum_{i\neq 0} P_{0,0}(S_n = i)P_{0,i}(\tilde{H}_0 < \infty) \leq M\delta_1^n.$$

Indeed, by conditioning on the last time the walk meets 0 in the period $[0,n]$, we obtain

(A.15) $$\sum_{i\neq 0} P_{0,0}(S_n = i)P_{0,i}(\tilde{H}_0 < \infty)$$
$$\leq \sum_{0\leq k<n} P_{0,0}(S_k = 0)P_{0,0}(n-k < \tilde{H}_0 < \infty)$$
$$\leq 1/\delta_1^n \sum_{0\leq k<n} P_{0,0}(S_k = 0)\delta_1^k.$$

Now, note that, for $k \geq 1$,

$$P_{0,0}(S_k = 0) = \sum_{i=1}^{k-1} P_{0,0}(\tilde{H}_0 = i) P_{0,0}(S_{k-i} = 0) + P_{0,0}(\tilde{H}_0 = k),$$

so that, by (A.7),

$$M := \sum_{k \geq 0} P_{0,0}(S_k = 0) \delta_1^k = \frac{E[\delta_1^{\tilde{H}_0} \mathbb{1}\{\tilde{H}_0 < \infty\}]}{1 - E[\delta_1^{\tilde{H}_0} \mathbb{1}\{\tilde{H}_0 < \infty\}]} < \infty.$$

It is now easy to see how estimate (v) follows.

The last property is seen by first applying the Cauchy–Schwarz inequality as, for instance, for the discrete version of (ii):

$$\left( \sum_i \varepsilon_i P_{0,i}(\tilde{H}_0 = n) \right)^2 \leq \left( \sum_i \varepsilon_i^2 \right) \sum_i P_{0,i}(\tilde{H}_0 = n)^2 \stackrel{n \to \infty}{\longrightarrow} 0 \qquad \text{for } d \geq 3. \quad \square$$

C.M.I.  
Université de Provence  
39 Rue Joliot–Curie  
F-13453 Marseille cedex 13  
France  
E-mail: asselah@cmi.univ-mrs.fr

IME-USP  
P.B. 66281  
05315-970  
São Paulo, SP  
Brazil  
E-mail: pablo@ime.usp.br